\documentclass[11pt, a4paper]{article}

\usepackage{amsmath}
\usepackage{amsthm}
\usepackage{amssymb}
\usepackage{mathrsfs}
\usepackage{dsfont}
\usepackage{color}
\usepackage{mathabx}

\newcommand{\tab}{\hspace*{2em}}

\begin{document}

\title{Ideal Theory in Rings\\
	(Idealtheorie in Ringbereichen)}
\author{Emmy Noether\\
	Translated by Daniel Berlyne}
\maketitle

{\LARGE \textbf{Contents.}} \\

\begin{tabular}{ll}
& Introduction. \\
\S 1. & Ring, ideal, finiteness condition. \\
\S 2. & Representation of an ideal as the least common multiple of finitely \\ & many irreducible ideals. \\
\S 3. & Equality of the number of components in two different \\ & decompositions into irreducible ideals. \\
\S 4. & Primary ideals. Uniqueness of the prime ideals belonging to two \\ & different 
decompositions into irreducible ideals. \\
\S 5. & Representation of an ideal as the least common multiple of maximal \\ & primary ideals. 
Uniqueness of the associated prime ideals. \\
\S 6. & Unique representation of an ideal as the least common multiple of \\ & relatively prime  
irreducible ideals. \\
\S 7. & Uniqueness of the isolated ideals. \\
\S 8. & Unique representation of an ideal as the product of coprime \\ & irreducible ideals. \\
\S 9. & Development of the study of modules. Equality of the number of \\ & components in 
decompositions into irreducible modules. \\
\S 10. & Special case of the polynomial ring. \\
\S 11. & Examples from number theory and the theory of differential \\ & expressions. \\
\S 12. & Example from elementary divisor theory. \\
& Translator's notes. \\
& Acknowledgements. \\
\end{tabular}
\newpage

\emph{As this paper was originally written in the early twentieth century, there are a number of mathematical terms used that do not have an exact modern equivalent, and as such may be ambiguous in meaning. Such terms are underlined as they appear in the text, and an explanation is given at the end of the paper to clarify their meanings.}

\section*{Introduction.}
\tab This paper aims to \emph{convert the decomposition theorems for the integers or the decomposition of ideals in algebraic number fields into theorems for ideals in arbitrary integral domains (and rings in general)}. To understand this correspondence, we consider the decomposition theorems for the integers in a form somewhat different from the commonly given formulation. \\
\tab In the equation

\[ a = p_{1}^{\varrho_{1}}p_{2}^{\varrho_{2}} \dots p_{\sigma}^{\varrho_{\sigma}} = q_{1}q_{2} \dots q_{\sigma} \]
take the prime powers $q_{i}$ to be components of the decomposition with the following characteristic properties: \\
\tab 1. They are \underline{\emph{pairwise coprime}}; and no $q$ can be written as a product of pairwise coprime numbers, so in this sense irreducibility holds. \\ 
\tab 2. Each two components $q_{i}$ and $q_{k}$ are \underline{\emph{relatively prime}}; that is to say, if $bq_{i}$ is divisible by $q_{k}$, then $b$ is divisible by $q_{k}$. Irreducibility also holds in this sense. \\ 
\tab 3. Every $q$ is \emph{primary}; that is to say, if a product $bc$ is divisible by $q$, but $b$ is not divisible by $q$, then a power\footnote{If this power is always the first, then as is well-known it concerns a prime number.} of $c$ is divisible by $q$. The representation furthermore consists of \emph{maximal primary components}, since the product of two different $q$ is no longer primary. The $q$ are also irreducible in relation to the decomposition into maximal primary components. \\
\tab 4. Each $q$ is \emph{irreducible} in the sense that it cannot be written as the least common multiple of two proper divisors.

The connection between these primary numbers $q$ and the prime numbers $p$ is that for every $q$ there is one and (disregarding sign) only one $p$ that is a divisor of $q$ and a power of which is divisible by $q$: the associated prime number. If $p^{\varrho}$ is the lowest such power \---- $\varrho$ being the exponent from $q$ \---- then in particular, $p^{\varrho}$ is equal to $q$ here. The \emph{uniqueness theorem} can now be stated as follows: \\
\tab \emph{Given two different decompositions of an integer into irreducible maximal primary components $q$, each decomposition has the same number of components, the same associated prime numbers (up to sign) and the same exponents. Because $p^{\varrho} = q$, it also follows that the $q$ themselves are the same (up to sign).}

As is well-known, the uncertainty brought about by the sign is eradicated if instead of numbers, the ideals derived from them (all numbers divisible by $a$) are considered; the formulation then holds exactly for the unique decomposition of ideals of (finite) algebraic number fields into powers of prime ideals.

In the following (\S 1) a general ring is considered, which must only satisfy the \emph{finiteness condition} that every ideal of the domain has a finite ideal basis. Without such a finiteness condition, irreducible and prime ideals need not exist, as shown by the ring of \emph{all} algebraic integers, in which there is no decomposition into prime ideals.

It is shown that \---- corresponding to the four characteristic properties of the component $q$ \---- in general \emph{four separate decompositions exist, which follow successively from each other through subdivision}. Thereby the decomposition into coprime irreducible ideals behaves as a factorisation, and the remaining three decompositions behave as a reduced (\S 2) representation as the least common multiple. Also, the connection between a primary ideal \---- the irreducible ideals are also primary \---- and the corresponding prime ideal is preserved: every primary ideal $\mathfrak{Q}$ uniquely determines a corresponding prime ideal $\mathfrak{P}$ which is a divisor of $\mathfrak{Q}$ and a power of which is divisible by $\mathfrak{Q}$.  If $\mathfrak{P}^{\varrho}$ is the lowest such power \---- $\varrho$ being the exponent of $\mathfrak{Q}$ \---- then $\mathfrak{P}^{\varrho}$ does not need to coincide with $\mathfrak{Q}$ here, however. The \emph{uniqueness theorem} can now be expressed as follows:

\emph{The decompositions $1$ and $2$ are unique; given two different decompositions $3$ or $4$, the number of components and the corresponding prime ideals are the same;\footnote{Moreover, the exponents are presumably also the same, and more generally the corresponding components are isomorphic.} the isolated ideals} (\S 7) \emph{occurring in the components are uniquely determined.}

For the proofs of the decomposition theorems, the notable \textquotedblleft Theorem of the Finite Chain" for finite modules, first by Dedekind, is followed using the finiteness condition, and from this we deduce representation $4$ of each ideal as the least common multiple of finitely many irreducible ideals. By revising the statement of the reducibility of a component, the fundamental \emph{uniqueness theorem for decomposition $4$ into irreducible ideals} is produced. By concluding that each of the remaining decompositions are given by finitely many components, the uniqueness theorems for these emerge as a result of uniqueness theorem $4$.

Finally it is shown (\S 9) that the representation through finitely many \emph{irreducible} components also holds under weaker requirements; the commutativity of the ring is not necessary, and it suffices to consider a module in relation to the ring instead of an ideal. In this more general case the \emph{equality of the number of components} for two different decompositions still holds, while the notions of prime and primary are restricted to commutativity and the concept of an ideal; in contrast, the notion of coprime for ideals is retained in non-commutative rings.

The simplest ring for which the four separate decompositions actually occur is the ring of all polynomials in $n$ variables with arbitrary complex coefficients. The particular decompositions can be deduced here to be irrational due to the properties of algebraic figures, and the uniqueness theorem for the corresponding prime ideals is equivalent to the Fundamental Theorem of Elimination Theory regarding the unique decomposability of algebraic figures into irreducible elements.
Further examples are given by all finite integral domains of polynomials (\S 10). In fact, the simple ring of all even numbers, or more generally all numbers divisible by a given number, is also an example for the separate decompositions (\S 11). An example of ideal theory in non-commutative rings is provided by elementary divisor theory (\S 12), where unique decomposition into irreducible ideals, or classes, holds. These irreducible classes characterise completely the irreducible parts of elementary divisors, and in rings where the usual elementary divisor theory breaks down can perhaps be considered as their equivalent.

In the available literature the following are to be noted: the decomposition into \emph{maximal primary ideals} is given by Lasker for the polynomial ring with arbitrary complex or integer coefficients, and taken further by Macaulay at particular points.\footnote{E. Lasker, Zur Theorie der Moduln und Ideale. Math. Ann. 60 (1905), p20, Theorems VII and XIII. \---- F. S. Macaulay, On the Resolution of a given Modular System into Primary Systems including some Properties of Hilbert Numbers. Math. Ann. 74 (1913), p66.} Both concern themselves with elimination theory, therefore using the fact that a polynomial can be expressed uniquely as the product of irreducible polynomials. In fact, the decomposition theorems for ideals are independent of this hypothesis, as ideal theory in algebraic number fields allows one to suppose and as this paper shows. The primary ideal is also defined by Lasker and Macaulay using concepts from elimination theory.

The decomposition into \emph{irreducible} ideals and into \emph{relatively prime irreducible} ideals appears also not to be remarked upon for the polynomial ring in the available literature; only a remark by Macaulay on the uniqueness of the isolated primary ideals can be found.

The decomposition into \emph{coprime irreducible} ideals is given by Schmeidler\footnote{W. Schmeidler, \"{U}ber Moduln und Gruppen hyperkomplexer Gr\"{o}\ss en. Math. Zeitschr. 3 (1919), p29.} for the polynomial ring, using elimination theory for the proof of finiteness. However, here the uniqueness theorem is stated only for classes of ideals, not for the ideals themselves. This last uniqueness theorem can be found in a joint paper,\footnote{E. Noether - W. Schmeidler, Moduln in nichtkommutativen Bereichen, insbesondere aus Differential- und Differenzenausdr\"{u}cken. Math. Zeitschr. 8 (1920), p1.} where it refers to ideals in non-commutative polynomial rings. Only the finite ideal basis is made use of here, so theorems and methods for general rings remain to be addressed, becoming more of a problem through this paper in terms of equality in size (\S 11). The present researches give a strong generalisation and further development of the underlying concepts of both of these works. The basis of both works is the transition from the expression as the least common multiple to an additive decomposition of the system of residue classes. Here the least common multiple is kept for the sake of a simpler representation; the additive decomposition then corresponds to the conversion of the idea of reducibility into a property of the complement (\S 3). Therefore all given theorems are left to be reflected upon and understood in the form of additive decomposition theorems for the system of residue classes and known subsystems, which is essentially equivalent to the reflections of the joint paper. This system of residue classes forms a ring of the same generality as originally laid down; that is, every ring can be regarded as a system of residue classes of the ideal which corresponds to the collection of all identity relations between the elements of the ring; or also a subsystem of these relations, by assuming the remaining relations are also satisfied in the ring.

This remark also gives the classification in Fraenkel's papers.\footnote{A. Fraenkel, \"{U}ber die Teiler der Null und die Zerlegung von Ringen. J. f. M. 145 (1914), p139. \"{U}ber gewisse Teilbereiche und Erweiterungen von Ringen. Professorial dissertation, Leipzig, Teubner, 1916. \"{U}ber einfache Erweiterungen zerlegbarer Ringe. J. f. M. 151 (1920), p121.} Fraenkel considers additive decompositions of rings that depend on such restrictive conditions (existence of regular elements, division by these, decomposability requirement), which coincide with the four decompositions for the corresponding ideal. Because of this concurrence, its finiteness condition also means that the ideal only has finitely many proper divisors \---- apart from some exceptional cases \---- a restriction no stricter than ours. Fraenkel's starting point is different, dependent on the essentially algebraic goals of his work; through algebraic extension he achieved more general rings with fewer restricting conditions.

\section*{\S 1. Ring, ideal, finiteness requirement.}

\tab \textbf{1.} Let $\Sigma$ be a (commutative) \emph{ring} in an abstract definition;\footnote{The definition is taken from Fraenkel's professorial dissertation, with omission of the more restrictive requirements 6, I and II; instead the commutative law of addition must be incorporated. It therefore concerns the laws defining a field, with the omission of the multiplicative inverse.} that is to say, $\Sigma$ consists of a system of elements $a, b, c, \dots, f, g, h, \dots$ in which a relation satisfying the usual requirements is defined as \emph{equality}, and in which each two ring elements $a$ and $b$ combine uniquely through two operations, \emph{addition} and \emph{multiplication}, to give a third element, given by the sum $a + b$ and the product $a \cdot b$. The ring and the otherwise entirely arbitrary operations must satisfy the following laws: \\
\tab 1. \emph{The associative law of addition:} $(a+b)+c = a+(b+c).$ \\
\tab 2. \emph{The commutative law of addition:} $ a+b = b+a.$ \\
\tab 3. \emph{The associative law of multiplication:} $(a \cdot b) \cdot c = a \cdot (b \cdot c).$ \\
\tab 4. \emph{The commutative law of multiplication:} $a \cdot b = b \cdot a.$ \\
\tab 5. \emph{The distributive law:} $a \cdot (b+c) = a \cdot b + a \cdot c.$ \\
\tab 6. \emph{The law of unrestricted and unique subtraction.} \\
In $\Sigma$ there is a single element $x$ which satisfies the equation $a+x = b$. (Written $x = b-a$.)

The existence of the zero element follows from these properties; however, a ring is not required to possess a unit, and the product of two elements can be zero without either of the factors being zero. Rings for which a product being zero implies one of the factors is zero, and which in addition possess a unit, are called \emph{integral domains}. For the finite sum $a+a+ \dots +a$ we use the usual abbreviated notation $na$, where the integers $n$ are solely considered as shortened notation, not as ring elements, and are defined recursively through $a = 1 \cdot a$, $na+a = (n+1)a$.

\textbf{2.} Let an \emph{ideal} $\mathfrak{M}$\footnote{Ideals are denoted by capital letters. The use of $\mathfrak{M}$ brings to mind the example of the ideals composed of polynomials usually called \textquotedblleft modules". Incidentally, \S\S 1-3 use only the properties of modules and not the properties of ideals; compare \S 9 as well.} in $\Sigma$ be understood to be a system of elements of $\Sigma$ such that the following two conditions are satisfied: \\
\tab 1. \emph{If} $\mathfrak{M}$ \emph{contains} $f$\emph{, then} $\mathfrak{M}$ \emph{also contains} $a \cdot f$\emph{, where} $a$ \emph{is an arbitrary element of} $\Sigma$. \\
\tab 2. \emph{If} $\mathfrak{M}$ \emph{contains} $f$ \emph{and} $g$\emph{, then} $\mathfrak{M}$ \emph{also contains the difference} $f-g$; so if $\mathfrak{M}$ contains $f$, then $\mathfrak{M}$ also contains $nf$ for all integers $n$.

If $f$ is an element of $\mathfrak{M}$, then we write $f \equiv 0\ (\mathfrak{M})$, as is usual; and we say that $f$ \emph{is divisible by} $\mathfrak{M}$. If every element of $\mathfrak{N}$ is equal to some element of $\mathfrak{M}$, and so divisible by $\mathfrak{M}$, then we say $\mathfrak{N}$ \emph{is divsible by} $\mathfrak{M}$, written $\mathfrak{N} \equiv 0\ (\mathfrak{M})$. $\mathfrak{M}$ is called a \emph{proper divisor} of $\mathfrak{N}$ if it contains elements not in $\mathfrak{N}$, and so is not conversely divisible by $\mathfrak{N}$. If $\mathfrak{N} \equiv 0\ (\mathfrak{M})$ and $\mathfrak{M} \equiv 0\ (\mathfrak{N})$, then $\mathfrak{N} = \mathfrak{M}$.

The remaining familiar notions also remain valid, word for word. By the \emph{greatest common divisor} $\mathfrak{D} = (\mathfrak{A}, \mathfrak{B})$ of two ideals $\mathfrak{A}$ and $\mathfrak{B}$, we understand this to be all elements which are expressible in the form $a+b$, where $a$ is an element of $\mathfrak{A}$ and $b$ is an element of $\mathfrak{B}$; $\mathfrak{D}$ is also itself an ideal. Likewise, the greatest common divisor $\mathfrak{D} = (\mathfrak{A}_{1}, \mathfrak{A}_{2}, \dots, \mathfrak{A}_{\nu}, \dots )$ of infinitely many ideals is defined as all elements $d$ which are expressible as the sum of the elements of each of finitely many ideals: $d = a_{i_{1}} + a_{i_{2}} + \dots + a_{i_{n}}$; here $\mathfrak{D}$ is again itself an ideal.

Should the ideal $\mathfrak{M}$ contain in particular a finite number of elements $f_{1}, f_{2}, \dots, f_{\varrho}$ such that
\begin{align*}
\mathfrak{M} = (f_{1}, \dots, f_{\varrho})\text{;}
\end{align*}
that is,
\begin{align*}
f = a_{1}f_{1}+ \dots +a_{\varrho}f_{\varrho}+n_{1}f_{1}+ \dots n_{\varrho}f_{\varrho}\text{ for all } f \equiv 0\ (\mathfrak{M})\text{,}
\end{align*}
where the $a_{i}$ are elements of the \emph{ring}, the $n_{i}$ are integers, and so $\mathfrak{M}$ forms a \emph{finite ideal}; $f_{1}, \dots, f_{\varrho}$ forms an \emph{ideal basis}.

In the following we now consider solely \emph{rings} $\Sigma$ \emph{which satisfy the finiteness condition: every ideal in} $\Sigma$ \emph{is finite, and so has an ideal basis}.

\textbf{3.} The following underlying ideas all follow directly from the finiteness condition:

\textbf{Theorem I} (\emph{Theorem of the Finite Chain}):\footnote{Initially stated for modules by Dedekind: Zahlentheorie, Suppl. XI, \S 172, Theorem VIII (4th condition); our proof and the term "chain" is taken from there. For ideals of polynomials by Lasker, loc. cit. p56 (lemma). In both cases the theorem finds only specific applications, however. Our applications depend without exception on the axiom of choice.} \emph{If} $\mathfrak{M}, \mathfrak{M}_{1}, \mathfrak{M}_{2}, \dots, \mathfrak{M}_{\nu}, \dots$ \emph{is a countably infinite system of ideals in} $\Sigma$ \emph{in which each ideal is divisible by the following one, then all ideals after a finite index n are identical;} $\mathfrak{M}_{n} = \mathfrak{M}_{n+1} = \dots$. \emph{In other words: If} $\mathfrak{M}, \mathfrak{M}_{1}, \mathfrak{M}_{2}, \dots, \mathfrak{M}_{\nu}, \dots$ \emph{gives a simply ordered chain of ideals such that each ideal is a proper divisor of its immediate predecessor, then the chain terminates in a finite number of steps.} \\

In particular, let $\mathfrak{D} = (\mathfrak{M}_{1}, \mathfrak{M}_{2}, \dots, \mathfrak{M}_{\nu}, \dots)$ be the greatest common divisor of the system, and let $f_{1} \dots f_{k}$ be an always existing basis of $\mathfrak{D}$ resulting from the finiteness condition. Then it follows from the divisibility assumption that every element of $\mathfrak{D}$ is equal to an element of an ideal in the chain; 
then it follows from
\begin{align*}
f = g+h,\quad g \equiv 0\ (\mathfrak{M}_{r}),\quad h \equiv 0\ (\mathfrak{M}_{s}),\quad (r \leq s)
\end{align*}
that $g \equiv 0\ (\mathfrak{M}_{s})$ and so $f \equiv 0\ (\mathfrak{M}_{s})$. The corresponding statement holds if $f$ is the sum of several components. There is therefore also a finite index $n$ such that
\begin{align*}
f_{1} \equiv 0\ (\mathfrak{M}_{n});\ \dots;\ f_{k} \equiv 0\ (\mathfrak{M}_{n});\ \mathfrak{D} = (f_{1}, \dots, f_{k}) \equiv 0\ (\mathfrak{M}_{n}).
\end{align*}
Because conversely $\mathfrak{M}_{n} \equiv 0\ (\mathfrak{D})$, it follows that $\mathfrak{M}_{n} = \mathfrak{D}$; and because furthermore
\begin{align*}
\mathfrak{M}_{n+i} \equiv 0\ (\mathfrak{D}); \quad \mathfrak{D} = \mathfrak{M}_{n} \equiv 0\ (\mathfrak{M}_{n+i}),
\end{align*}
it also follows that $\mathfrak{M}_{n+i} = \mathfrak{D} = \mathfrak{M}_{n}$ for all $i$, whereupon the theorem is proved. \\

Note that conversely the existence of the ideal basis follows from this theorem, so \emph{the finiteness condition could also have been stated in this basis-free form}.

\section*{\S 2. Representation of an ideal as the least common multiple of finitely many irreducible ideals.}

\tab Let the \emph{least common multiple} $[\mathfrak{B}_{1}, \mathfrak{B}_{2}, \dots, \mathfrak{B}_{k}]$ of the ideals $\mathfrak{B}_{1}, \mathfrak{B}_{2},$ $\dots, \mathfrak{B}_{k}$ be defined as usual as the collection of elements which are divisible by each of $\mathfrak{B}_{1}, \mathfrak{B}_{2}, \dots, \mathfrak{B}_{k}$, written as:
\begin{align*}
f \equiv 0\ (\mathfrak{B}_{i}),\ (i = 1, 2, \dots, k)\ \text{implies } f \equiv 0\ ([\mathfrak{B}_{1}, \mathfrak{B}_{2}, \dots, \mathfrak{B}_{k}])
\end{align*}
and vice versa. The least common multiple is itself an ideal; we call the ideals $\mathfrak{B}_{i}$ the \emph{components of the decomposition}. \\

\textbf{Definition I}. \emph{A representation} $\mathfrak{M} = [\mathfrak{B}_{1}, \dots, \mathfrak{B}_{k}]$ \emph{is called a reduced representation if no} $\mathfrak{B}_{i}$ \emph{appears in the least common multiple} $\mathfrak{A}_{i}$ \emph{of the remaining ideals, and if no} $\mathfrak{B}_{i}$ \emph{can be replaced with a proper divisor}.\footnote{An example of a non-reduced representation is: 
\begin{align*}
(x^{2},\ xy) = [(x),\ (x^{2},\ xy,\ y^{\lambda})]
\end{align*}
for all exponents $\lambda \geq 2$; the case $\lambda = 1$, corresponding to the representation $[(x),\ (x^{2},\ y)]$, gives a reduced representation. (K. Hentzelt, who was killed in the war, gave this representation to me as the simplest example of a non-unique decomposition into primary ideals.)} \emph{If the conditions are only satisfied for the ideal} $\mathfrak{B}_{i}$\emph{, then the representation is called reduced with respect to} $\mathfrak{B}_{i}$. \emph{The least common multiple} $\mathfrak{A}_{i} = [\mathfrak{B}_{1}, \dots, \mathfrak{B}_{i-1}, \mathfrak{B}_{i+1}, \dots, \mathfrak{B}_{k}]$ \emph{is called the complement of} $\mathfrak{B}_{i}$. \emph{Representations in which only the first condition is satisfied are called shortest representations}. \\

It suffices now, when considering the representation of an ideal as a lowest common multiple, to restrict ourselves to reduced representations, due to the following lemma: \\

\textbf{Lemma I}. \emph{Every representation of an ideal as the least common multiple of finitely many ideals can be replaced in at least one way by a reduced representation; in particular, one such representation can be obtained through successive decomposition.} \\

Let $\mathfrak{M} = [\mathfrak{B}_{1}^{*}, \dots, \mathfrak{B}_{l}^{*}]$ be an arbitrary representation of $\mathfrak{M}$. We can then omit those $\mathfrak{B}_{i}^{*}$ which go into the least common multiple of the remaining ideals. Because the remaining ideals still give $\mathfrak{M}$, the resulting representation
\begin{align*}
\mathfrak{M} = [\mathfrak{B}_{1}, \dots, \mathfrak{B}_{k}] = [\mathfrak{A}_{i}, \mathfrak{B}_{i}]
\end{align*}
satisfies the first condition, and so this is a shortest representation; and this condition remains satisfied if some $\mathfrak{B}_{i}$ is replaced with a proper divisor. But the second condition is always satisfiable, by the Theorem of the Finite Chain (Theorem I). For suppose
\begin{align*}
\mathfrak{M} = [\mathfrak{A}_{i}, \mathfrak{B}_{i}] = [\mathfrak{A}_{i}, \mathfrak{B}_{i}^{(1)}] = \dots = [\mathfrak{A}_{i}, \mathfrak{B}_{i}^{(\nu)}], \dots,
\end{align*}
where each $\mathfrak{B}_{i}^{(\nu)}$ is a proper divisor of its immediate predecessor, so that the chain $\mathfrak{B}_{i}, \mathfrak{B}_{i}^{(1)}, \dots, \mathfrak{B}_{i}^{(\nu)}, \dots$ must terminate in finitely many steps; in the representation $\mathfrak{M} = [\mathfrak{A}_{i}, \mathfrak{B}_{i}^{(n)}]$, $\mathfrak{B}_{i}^{(n)}$ can therefore not be replaced with a proper divisor, and this holds a fortiori if $\mathfrak{A}_{i}$ is replaced with a proper divisor. Therefore if the algorithm is applied to each $\mathfrak{B}_{i}$ by at each stage using the complement together with the already reduced $\mathfrak{B}$, then a reduced representation is formed.\footnote{The previous example shows that the given representation does not uniquely define such a reduced representation. For $(x^{2},\ xy) = [(x),\ (x^{2},\ xy,\ y^{\lambda})]$, where $\lambda \geq 2$, both $[(x),\ (x^{2},\ y)]$ and $[(x),\ (x^{2},\ \mu x + y)]$ for arbitrary $\mu$ are reduced representations.}

In order to obtain such a representation successively, it must be shown that it follows from the individual reduced representations
\begin{align*}
\mathfrak{M} = [\mathfrak{B}_{1}, \mathfrak{C}_{1}], \mathfrak{C}_{1} = [\mathfrak{B}_{2}, \mathfrak{C}_{2}],\ \dots,\ \mathfrak{C}_{l-1} = [\mathfrak{B}_{l}, \mathfrak{C}_{l}]
\end{align*}
that the representation $\mathfrak{M} = [\mathfrak{B}_{1}, \dots, \mathfrak{B}_{l}, \mathfrak{C}_{l}]$ resulting from these is also reduced. For this it is sufficient to show that if the representations $\mathfrak{M} = [\mathfrak{B}_{1}, \dots, \mathfrak{B}_{\varrho}, \mathfrak{C}]$ and $\mathfrak{C} = [\mathfrak{C}_{1}, \mathfrak{C}_{2}]$ are reduced, then $\mathfrak{M} = [\mathfrak{B}_{1}, \dots, \mathfrak{B}_{\varrho}, \mathfrak{C}_{1}, \mathfrak{C}_{2}]$ is also reduced. Indeed, by assumption no $\mathfrak{B}$ appears in its complement; were this the case for a $\mathfrak{C}_{i}$, then in the first representation, $\mathfrak{C}$ would be replaceable with a proper divisor, contrary to the assumption, because due to the second reduced representation, $\mathfrak{C}_{1}$ and $\mathfrak{C}_{2}$ are proper divisors of $\mathfrak{C}$; the representation is therefore a shortest representation. Furthermore, by assumption no $\mathfrak{B}$ can be replaced with a proper divisor; were this the case for a $\mathfrak{C}_{i}$, then this would correspond to the substitution of $\mathfrak{C}$ with a proper divisor, in contradiction with the assumption, because the representation for $\mathfrak{C}$ is reduced. \emph{The lemma is thus proved}. \\

\textbf{Definition II.} \emph{An ideal} $\mathfrak{M}$ \emph{is called reducible if it can be expressed as the least common multiple of two proper divisors; otherwise} $\mathfrak{M}$ \emph{is called irreducible.} \\

We now prove the following via the Theorem of the Finite Chain (Theorem I) and using the reduced representation: \\ 

\textbf{Theorem II.} \emph{Every ideal can be expressed as the least common multiple of finitely many irreducible ideals}.\footnote{The previous example shows that such a representation is in general not unique: $(x^{2}, xy) = [(x), (x^{2}, \mu x + y)]$. Both components are irreducible for arbitrary $\mu$. All divisors of $(x)$ are of the form $(x, g(y))$, where $g(y)$ denotes a polynomial in $y$; therefore if the least common multiple of two divisors also has this form, then it is a proper divisor of $(x)$. $(x^{2}, \mu x + y)$ has only \emph{one} divisor $(x, y)$, and so is necessarily also irreducible.} \\

An arbitrary ideal $\mathfrak{M}$ is either irreducible, in which case $\mathfrak{M} = [\mathfrak{M}]$ is a representation in the form required by Theorem II, or it is $\mathfrak{M} = [\mathfrak{B}_{1}, \mathfrak{C}_{1}]$, where $\mathfrak{B}_{1}, \mathfrak{C}_{1}$ are proper divisors of $\mathfrak{M}$, and by Lemma I the representation can be assumed to be reduced. The same choice is true of $\mathfrak{C}_{1}$; either it is irreducible or it has a reduced representation $\mathfrak{C}_{1} = [\mathfrak{B}_{2}, \mathfrak{C}_{2}]$.

Continuing in this way, the following series of reduced representations is obtained:
\begin{align}
\mathfrak{M} = [\mathfrak{B}_{1}, \mathfrak{C}_{1}];\ \mathfrak{C}_{1} = [\mathfrak{B}_{2}, \mathfrak{C}_{2}];\ \dots ;\ \mathfrak{C}_{\nu - 1} = [\mathfrak{B}_{\nu}, \mathfrak{C}_{\nu}];\ \dots \quad \text{.}
\end{align}

In the chain $\mathfrak{C}_{1}, \mathfrak{C}_{2}, \dots , \mathfrak{C}_{\nu}, \dots$, each $\mathfrak{C}_{i}$ is a proper divisor of its immediate predecessor, and therefore the chain terminates in finitely many steps; there is an index $n$ such that $\mathfrak{C}_{n}$ is irreducible. Furthermore, by Lemma I the representation $\mathfrak{M} = [\mathfrak{B}_{1}, \dots, \mathfrak{B}_{n}, \mathfrak{C}_{n}]$ is reduced; $\mathfrak{C}_{n}$ can therefore not go into its complement $\mathfrak{A}_{n}$, and $\mathfrak{C}_{n}$ cannot be replaced by any proper divisor in the representation $\mathfrak{M} = [\mathfrak{A}_{n}, \mathfrak{C}_{n}]$. Replacing $\mathfrak{A}_{n}$ with a proper divisor\footnote{Indeed, the representation is also reduced with respect to $\mathfrak{A}_{n}$, as will be shown in \S 3 (Lemma IV) as the converse of Lemma I.} if necessary so that the representation is reduced, it has therefore been shown that \emph{every reducible ideal admits a reduced representation as the least common multiple of an irreducible ideal and an ideal complementary to it}. All $\mathfrak{B}_{i}$ in series (1) can therefore without loss of generality be assumed to be irreducible. Iteration of the above argument gives the existence of an irreducible $\mathfrak{C}_{n}$, which \emph{proves} \emph{Theorem} II.

\section*{\S 3. Equality of the number of components in two different decompositions into irreducible ideals.}

\tab To prove the equality of the number of components, it is first necessary to express the reducibility and irreducibility of an ideal through properties of its complement, via \\

\textbf{Theorem III}.\footnote{Theorem III corresponds to the transition from modules to quotient groups in the works of Schmeidler and Noether-Schmeidler (cf. the introduction). Here $\mathfrak{A}$ corresponds to the quotient group, and $\mathfrak{N}_{1}$ and $\mathfrak{N}_{2}$ to the subgroups into which the quotient group is decomposed.} \emph{Let the shortest representation} $\mathfrak{M} = [\mathfrak{A},\mathfrak{C}]$ \emph{be reduced with respect to} $\mathfrak{C}$. \emph{Then a necessary and sufficient condition for} $\mathfrak{C}$ \emph{to be reducible is the existence of two ideals} $\mathfrak{N}_{1}$ \emph{and} $\mathfrak{N}_{2}$ \emph{which are proper divisors of} $\mathfrak{M}$ \emph{such that}
\begin{align}
\mathfrak{N}_{1} \equiv 0\ (\mathfrak{A}); \quad \mathfrak{N}_{2} \equiv 0\ (\mathfrak{A}); \quad [\mathfrak{N}_{1}, \mathfrak{N}_{2}] = \mathfrak{M}.
\end{align}
From this also follows: \emph{If the conditions} (2) \emph{are satisfied, and} $\mathfrak{C}$ \emph{is irreducible, then at least one of the} $\mathfrak{N}_{i}$ \emph{is not a proper divisor of} $\mathfrak{M}$; $\mathfrak{N}_{i} = \mathfrak{M}$. \\

Let $\mathfrak{C} = [\mathfrak{C}_{1}, \mathfrak{C}_{2}]$, where $\mathfrak{C}_{1},\ \mathfrak{C}_{2}$ are proper divisors of $\mathfrak{C}$. Then it follows that
\begin{align*}
\mathfrak{M} = [\mathfrak{A}, \mathfrak{C}] = [\mathfrak{A}, \mathfrak{C}_{1}, \mathfrak{C}_{2}] = [[\mathfrak{A}, \mathfrak{C}_{1}], [\mathfrak{A}, \mathfrak{C}_{2}]].
\end{align*}
Here the ideals $[\mathfrak{A}, \mathfrak{C}_{i}]$ are proper divisors of $\mathfrak{M}$, because otherwise $[\mathfrak{A}, \mathfrak{C}]$ would not be reduced with respect to $\mathfrak{C}$. Because the divisibility by $\mathfrak{M}$ is also satisfied, \emph{condition} (2) \emph{is proved to be necessary}. (The representation (2) is not reduced, because a $[\mathfrak{A}, \mathfrak{C}_{i}]$ can be replaced with $\mathfrak{C}_{i}$.)

Conversely, now suppose (2) holds. We construct the ideals:
\begin{align*}
\mathfrak{C}_{1} = (\mathfrak{C}, \mathfrak{N}_{1}); \quad \mathfrak{C}_{2} = (\mathfrak{C}, \mathfrak{N}_{2}); \quad \mathfrak{C}^{*} = [\mathfrak{C}_{1}, \mathfrak{C}_{2}].
\end{align*}
Then $\mathfrak{C}$ is divisible by both $\mathfrak{C}_{1}$ and $\mathfrak{C}_{2}$, and therefore also by the least common multiple $\mathfrak{C}^{*}$. In order to show the divisibility of $\mathfrak{C}^{*}$ by $\mathfrak{C}$, let $f \equiv 0\ (\mathfrak{C}^{*})$; therefore $f \equiv 0\ (\mathfrak{C}_{1})$ and $f \equiv 0\ (\mathfrak{C}_{2})$, or also $f = c + n_{1}$ and $f = \bar{c} + n_{2}$, where $c, \bar{c}, n_{1}, n_{2}$ are elements of $\mathfrak{C}, \mathfrak{N}_{1}, \mathfrak{N}_{2}$, and so in particular $n_{i}$ is divisible by $\mathfrak{A}$. Therefore the difference
\begin{align*}
g = c - \bar{c} = n_{2} - n_{1}
\end{align*}
is divisible both by $\mathfrak{C}$ and by $\mathfrak{A}$, and so by $\mathfrak{M}$. Because $n_{1} = n_{2} + m$, it further holds that $n_{1}$ (likewise $n_{2}$) is divisible both by $\mathfrak{N}_{1}$ and by $\mathfrak{N}_{2}$, and so by $\mathfrak{M}$. Therefore
\begin{align*}
f = c + m; \quad f \equiv 0\ (\mathfrak{C}); \quad \mathfrak{C}^{*} = \mathfrak{C}.
\end{align*}
Here $\mathfrak{C}_{1}$ and $\mathfrak{C}_{2}$ are proper divisors of $\mathfrak{C}$, since if $\mathfrak{C}_{i} = (\mathfrak{C}, \mathfrak{N}_{i}) = \mathfrak{C}$, then $\mathfrak{N}_{i}$ would be divisible by $\mathfrak{C}$, and so because of the divisibility by $\mathfrak{A}$, $\mathfrak{N}_{i}$ would be equal to $\mathfrak{M}$, in contradiction with the assumption. Hence $\mathfrak{C} = [\mathfrak{C}_{1}, \mathfrak{C}_{2}]$ is identified as reducible; \emph{Theorem} III \emph{is proved}.

Note that almost identical reasoning also shows the following: \\

\textbf{Lemma II.} \emph{If the ideal} $\mathfrak{C}$ \emph{in a shortest representation} $\mathfrak{M} = [\mathfrak{A}, \mathfrak{C}]$ \emph{can be replaced with a proper divisor, then} $\mathfrak{C}$ \emph{is reducible}. \\

Let
\begin{align*}
\mathfrak{M} = [\mathfrak{A}, \mathfrak{C}] = [\mathfrak{A}, \mathfrak{C}_{1}],
\end{align*}
and set
\begin{align*}
\mathfrak{C}^{*} = [\mathfrak{C}_{1}, (\mathfrak{A}, \mathfrak{C})].
\end{align*}
Then $\mathfrak{C}$ is again divisible by $\mathfrak{C}^{*}$. Furthermore, it follows from $f \equiv 0\ (\mathfrak{C}^{*})$ that
\begin{align*}
f = c_{1} = a + c.
\end{align*}
The difference $a = c_{1} - c$ is therefore divisible both by $\mathfrak{A}$ and by $\mathfrak{C}_{1}$, and therefore by $\mathfrak{M}$; therefore $c_{1} = c + m$; $f \equiv 0\ (\mathfrak{C})$; $\mathfrak{C}^{*} = \mathfrak{C}$.

Because both $\mathfrak{C}_{1}$ and $(\mathfrak{A}, \mathfrak{C})$ are proper divisors of $\mathfrak{C}$ by assumption, $\mathfrak{C} = \mathfrak{C}^{*}$ is thus identified as reducible.

An \emph{irreducible} $\mathfrak{C}$ can therefore not be replaced with a proper divisor.

Now let the following be \emph{two different shortest representations of} $\mathfrak{M}$ \emph{as the least common multiple of (finitely many) irreducible ideals}:
\begin{align*}
\mathfrak{M} = [\mathfrak{B}_{1}, \dots, \mathfrak{B}_{k}] = [\mathfrak{D}_{1}, \dots, \mathfrak{D}_{l}].
\end{align*}
These representations are both \emph{reduced} according to the remark following Lemma II. Now we shall first prove \\

\textbf{Lemma III.} \emph{For every complement} $\mathfrak{A}_{i} = [\mathfrak{B}_{1} \dots \mathfrak{B}_{i-1} \mathfrak{B}_{i+1} \dots \mathfrak{B}_{k}]$ \emph{there exists an ideal} $\mathfrak{D}_{j}$ \emph{such that} $\mathfrak{M} = [\mathfrak{A}_{i}, \mathfrak{D}_{j}]$. \\

Set $\mathfrak{M} = [\mathfrak{D}_{1}, \mathfrak{C}_{1}],\ \mathfrak{C}_{1} = [\mathfrak{D}_{2}, \mathfrak{C}_{12}]$ and so forth, so that
\begin{align*}
\mathfrak{M} = [\mathfrak{A}_{i}, \mathfrak{M}] = [\mathfrak{A}_{i}, \mathfrak{D}_{1}, \mathfrak{C}_{1}] = [[\mathfrak{A}_{i} \mathfrak{D}_{1}], [\mathfrak{A}_{i} \mathfrak{C}_{1}]].
\end{align*}
Here the conditions (2) from Theorem III are satisfied for $\mathfrak{N}_{1} = [\mathfrak{A}_{i}, \mathfrak{D}_{1}],\ \mathfrak{N}_{2} = [\mathfrak{A}_{i}, \mathfrak{C}_{1}]$, because $\mathfrak{M} = [\mathfrak{A}_{i}, \mathfrak{B}_{i}]$ is reduced with respect to $\mathfrak{B}_{i}$ and the representation is a shortest one. Because $\mathfrak{B}_{i}$ was assumed to be \emph{irreducible}, one $\mathfrak{N}_{i}$ must necessarily be equal to $\mathfrak{M}$.

For $\mathfrak{N}_{1} = \mathfrak{M}$ the lemma would be proved; for $\mathfrak{N}_{2} = \mathfrak{M}$ it follows respectively that $\mathfrak{M} = [[\mathfrak{A}_{i} \mathfrak{D}_{2}], [\mathfrak{A}_{i} \mathfrak{C}_{12}]]$, where, by the same result, one component must again be equal to $\mathfrak{M}$. Continuing in this way, it follows either that $\mathfrak{M} = [\mathfrak{A}_{i}, \mathfrak{D}_{j}]$, where $j < l$, or that $\mathfrak{M} = [\mathfrak{A}_{i}, \mathfrak{C}_{1 \dots l-1}]$; because $\mathfrak{C}_{1 \dots l-1} = \mathfrak{D}_{l}$, \emph{the lemma is thus proven}.

Now, as a result of this, we have \\

\textbf{Theorem IV.} \emph{For two different shortest representations of an ideal as the least common multiple of irreducible ideals, the number of components is the same.} \\

It follows from the lemma for the particular case $i = 1$:
\begin{align*}
\mathfrak{M} = [\mathfrak{A}_{1}, \mathfrak{B}_{1}] = [\mathfrak{A}_{1}, \mathfrak{D}_{j_1}] = [\mathfrak{D}_{j_1}, \mathfrak{B}_{2}, \dots, \mathfrak{B}_{k}].
\end{align*}
Now consider both of the decompositions
\begin{align*}
\mathfrak{M} = [\mathfrak{D}_{j_1}, \mathfrak{B}_{2}, \dots, \mathfrak{B}_{k}] = [\mathfrak{D}_{1}, \mathfrak{D}_{2}, \dots, \mathfrak{D}_{l}],
\end{align*}
and recall the earlier result with reference to the complement $\bar{\mathfrak{A}}_{2} =$ \\ $[\mathfrak{D}_{j_1}, \mathfrak{B}_{3}, \dots, \mathfrak{B}_{k}]$ of $\mathfrak{B}_{2}$. It then follows that
\begin{align*}
\mathfrak{M} = [\bar{\mathfrak{A}}_{2}, \mathfrak{B}_{2}] = [\bar{\mathfrak{A}}_{2}, \mathfrak{D}_{j_1}] = [\mathfrak{D}_{j_1}, \mathfrak{D}_{j_2}, \mathfrak{B}_{3}, \dots, \mathfrak{B}_{k}];
\end{align*}
and by extension of this procedure:
\begin{align*}
\mathfrak{M} = [\mathfrak{D}_{j_1}, \mathfrak{D}_{j_2}, \dots, \mathfrak{D}_{j_k}].
\end{align*}
Because now the representation $\mathfrak{M} = [\mathfrak{D}_{1}, \dots, \mathfrak{D}_{l}]$ is a shortest one by assumption, and so no $\mathfrak{D}$ can be omitted, the \emph{different} ones among the $\mathfrak{D}_{j_{i}}$ must exhaust all of the $\mathfrak{D}$; therefore it follows that $k \geq l$. Should we swap the $\mathfrak{B}$ with the $\mathfrak{D}$ throughout the lemma and the subsequent results, it then follows accordingly that $l \geq k$, and therefore that $k = l$, which proves the \emph{equality of the number of components}. As a result of this it follows that the ideals $\mathfrak{D}_{j_i}$ are all different from each other, because otherwise there would be fewer than $k$ components in a shortest representation using the $\mathfrak{D}_{j_i}$; the notation can therefore be chosen so that $\mathfrak{D}_{j_i} = \mathfrak{D}_{i}$. By the same reasoning, all intermediate representations $\mathfrak{M} = [\mathfrak{D}_{j_1}, \dots, \mathfrak{D}_{j_i}, \mathfrak{B}_{i+1}, \dots, \mathfrak{B}_{k}]$ are also shortest and so, by the remark following Lemma II, reduced.

The equality of the number of components leads us to a converse of Lemma I through

\textbf{Lemma IV.} \emph{If the components in a reduced representation are collected into groups and the least common multiples of them constructed, then the resulting representation is reduced. In other words: Given a reduced representation}
\begin{align*}
\mathfrak{M} = [\mathfrak{C}_{11}, \dots, \mathfrak{C}_{1 \mu_{1}}; \dots; \mathfrak{C}_{\sigma 1}, \dots, \mathfrak{C}_{\sigma \mu_{\sigma}}],
\end{align*}
\emph{it follows that} $\mathfrak{M} = [\mathfrak{N}_{1}, \dots, \mathfrak{N}_{\sigma}] = [\mathfrak{N}_{i}, \mathfrak{L}_{i}]$ \emph{is also reduced, where} $\mathfrak{N}_{i} = [\mathfrak{C}_{i1}, \dots, \mathfrak{C}_{i \mu_{i}}]$. \\

Firstly, note that $\mathfrak{N}_{i}$ cannot go into its complement $\mathfrak{L}_{i}$, since this is not the case for any of its divisors $\mathfrak{C}_{ij}$; the representation is therefore a shortest one. In order to show that $\mathfrak{N}_{i}$ cannot be replaced with any proper divisor, we split the $\mathfrak{C}$ into their irreducible ideals $\mathfrak{B}$,\footnote{By this we understand the $\mathfrak{B}$ to always be shortest, and therefore reduced, representations.} so that the (by Lemma I) reduced representations arise:
\begin{align*}
\mathfrak{M} = [\mathfrak{B}_{11}, \dots, \mathfrak{B}_{1 \lambda_{1}}; \dots; \mathfrak{B}_{\sigma 1}, \dots, \mathfrak{B}_{\sigma \lambda_{\sigma}}]; \quad \mathfrak{N}_{i} = [\mathfrak{B}_{i1}, \dots, \mathfrak{B}_{i \lambda_{i}}].
\end{align*}
Now let $\mathfrak{M} = [\mathfrak{N}_{i}^{*}, \mathfrak{L}_{i}]$ be reduced with respect to $\mathfrak{N}_{i}^{*}$, where $\mathfrak{N}_{i}^{*}$ is a \emph{proper} divisor of $\mathfrak{N}_{i}$. Then by Lemma II:
\begin{align*}
\mathfrak{N}_{i} = [\mathfrak{N}_{i}^{*}, (\mathfrak{N}_{i}, \mathfrak{L}_{i})];
\end{align*}
and this representation is necessarily reduced with respect to $\mathfrak{N}_{i}^{*}$, because otherwise $\mathfrak{N}_{i}^{*}$ can also be replaced with a proper divisor in $\mathfrak{M}$. Now also replace $(\mathfrak{N}_{i}, \mathfrak{L}_{i})$ with a proper divisor where appropriate, so that a reduced representation for $\mathfrak{N}_{i}$ is achieved. If both components of $\mathfrak{N}_{i}$ are now decomposed into irreducible ideals, then the number $\lambda_{i}$ of irreducible ideals in $\mathfrak{N}_{i}$ is composed additively of those of the components; the number of irreducible ideals corresponding to the proper divisor $\mathfrak{N}_{i}^{*}$ is therefore necessarily smaller than $\lambda_{i}$. Then, however, the decomposition of $\mathfrak{M} = [\mathfrak{N}_{i}^{*}, \mathfrak{L}]$ into irreducible ideals leads to fewer than $\sum_{i} \lambda_{i}$ ideals, \emph{in contradiction with the equality of number of components}. In the special case $\sigma = 2$, it also follows that the representation $\mathfrak{M} = [\mathfrak{N}_{i}, \mathfrak{L}_{i}]$ is reduced with respect to the complement $\mathfrak{L}_{i}$.

\section*{\S 4. Primary ideals. Uniqueness of the prime ideals belonging to two different decompositions into irreducible ideals.}

\tab The following concerns the connection between primary and irreducible ideals. \\

\textbf{Definition III.} \emph{An ideal} $\mathfrak{Q}$ \emph{is called primary if} $a \cdot b \equiv 0\ (\mathfrak{Q}),\ a \not\equiv 0\ (\mathfrak{Q})$ \emph{implies} $b^{x} \equiv 0\ (\mathfrak{Q})$, \emph{where the exponent} $x$ \emph{is a finite number}. \\

The definition can also be restated as follows: if a product $a \cdot b$ is divisible by $\mathfrak{Q}$, then either one factor is divisible by $\mathfrak{Q}$ or a power of the other is. \emph{If in particular} $x$ \emph{is always equal to} $1$, \emph{then the ideal is called a prime ideal}.

From the definition of a primary (respectively prime) ideal follows, by virtue of the existence of a basis, the definition using only products of ideals:\footnote{The product $\mathfrak{A} \cdot \mathfrak{B}$ of two ideals is understood to mean, as usual, the ideal consisting of the collection of elements $a \cdot b$ and their finite sums.} \\

\textbf{Definition IIIa.} \emph{An ideal} $\mathfrak{Q}$ \emph{is called primary if} $\mathfrak{A} \cdot \mathfrak{B} \equiv 0\ (\mathfrak{Q})$, $\mathfrak{A} \not\equiv 0\ (\mathfrak{Q})$ \emph{necessarily implies} $\mathfrak{B}^{\lambda} \equiv 0\ (\mathfrak{Q})$. \emph{If} $\lambda$ \emph{is always equal to} $1$, \emph{then the ideal is called a prime ideal. For a prime ideal} $\mathfrak{P}$, \emph{it therefore always follows from} $\mathfrak{A} \cdot \mathfrak{B} \equiv 0\ (\mathfrak{P})$, $\mathfrak{A} \not\equiv 0\ (\mathfrak{P})$ \emph{that} $\mathfrak{B} \equiv 0\ (\mathfrak{P})$. \\

Because Definition III is contained in IIIa for $\mathfrak{A} = (a)$, $\mathfrak{B} = (b)$ as a special case, every ideal which is primary by IIIa is also primary by III. Conversely, suppose $\mathfrak{Q}$ is primary by III, and let the assumption of IIIa be satisfied: $\mathfrak{A} \cdot \mathfrak{B} \equiv 0\ (\mathfrak{Q})$, so that it therefore follows either that $\mathfrak{A} \equiv 0\ (\mathfrak{Q})$, or alternatively that there is at least one element $a \equiv 0\ (\mathfrak{A})$ such that $a \cdot \mathfrak{B} \equiv 0\ (\mathfrak{Q})$ and $a \not\equiv 0\ (\mathfrak{Q})$. If now $b_{1}, \dots, b_{r}$ is an ideal basis of $\mathfrak{B}$, then by Definition III, since $a \cdot b_{i} \equiv 0\ (\mathfrak{Q})$, the following holds:
\begin{align*}
b_{1}^{x_{1}} \equiv 0\ (\mathfrak{Q});\ \dots;\ b_{r}^{x_{r}} \equiv 0\ (\mathfrak{Q}).
\end{align*}
Because $b = f_{1}b_{1} + \dots + f_{r}b_{r} + n_{1}b_{1} + \dots + n_{r}b_{r}$, for $\lambda = x_{1} + \dots + x_{r}$ the product of $\lambda$ elements $b$ is therefore divisible by $\mathfrak{Q}$, which proves the fulfilment of Definition IIIa for ideals primary by III. In particular, for prime ideals $\mathfrak{P}$ it follows from $a \cdot \mathfrak{B} \equiv 0\ (\mathfrak{P})$, and therefore also $a \cdot b \equiv 0\ (\mathfrak{P})$ for every $b \equiv 0\ (\mathfrak{B})$ and $a \not\equiv 0\ (\mathfrak{P})$, that $b \equiv 0\ (\mathfrak{P})$ and thus $\mathfrak{B} \equiv 0\ (\mathfrak{P})$. We have therefore shown the two definitions to be equivalent.

The connection between primary and prime ideals shall be established through the remark that the collection $\mathfrak{P}$ of all elements $p$ with the property that a power of $p$ is divisible by $\mathfrak{Q}$ forms a \emph{prime ideal}. It is immediately clear that $\mathfrak{P}$ is an ideal; because if the given property holds for $p_{1}$ and $p_{2}$, it also holds for $ap_{1}$ and $p_{1} - p_{2}$. Furthermore, by the inference used in Definition IIIa regarding the basis,
there exists a number $\lambda$ such that $\mathfrak{P}^{\lambda} \equiv 0\ (\mathfrak{Q})$.

Now let
\begin{align*}
a \cdot b \equiv 0\ (\mathfrak{P}); \quad a \not\equiv 0\ (\mathfrak{P}),
\end{align*}
so that it follows from the definition of $\mathfrak{P}$ that
\begin{align*}
a^{\lambda} \cdot b^{\lambda} \equiv 0\ (\mathfrak{Q}); \quad a^{\lambda} \not\equiv 0\ (\mathfrak{Q});
\end{align*}
therefore by the definition of $\mathfrak{Q}$:
\begin{align*}
b^{\lambda x} \equiv 0\ (\mathfrak{Q}) \quad \text{and hence} \quad b \equiv 0\ (\mathfrak{P}),
\end{align*}
which \emph{proves} $\mathfrak{P}$ \emph{to be a prime ideal}. $\mathfrak{P}$ is also defined as the greatest common divisor of all ideals $\mathfrak{B}$ with the property that a power of $\mathfrak{B}$ is divisible by $\mathfrak{Q}$. Every such $\mathfrak{B}$ is by definition divisible by $\mathfrak{P}$; thus the greatest common divisor $\mathfrak{D}$ of these $\mathfrak{B}$ is too. Conversely, $\mathfrak{P}$ is itself one of the ideals $\mathfrak{B}$, so is divisible by $\mathfrak{D}$, which proves that $\mathfrak{P} = \mathfrak{D}$. $\mathfrak{P}$ is therefore a prime ideal which is a divisor of $\mathfrak{Q}$ and a power of which is divisible by $\mathfrak{Q}$. It is uniquely defined by these properties, because it follows from
\begin{align*}
\mathfrak{Q} \equiv 0\ (\mathfrak{P}); \quad \mathfrak{P}^{\lambda} \equiv 0\ (\mathfrak{Q}); \quad \mathfrak{Q} \equiv 0\ (\bar{\mathfrak{P}}); \quad \bar{\mathfrak{P}}^{\mu} \equiv 0\ (\mathfrak{Q})
\end{align*}
that
\begin{align*}
\mathfrak{P}^{\lambda} \equiv 0\ (\bar{\mathfrak{P}}); \quad \bar{\mathfrak{P}}^{\mu} \equiv 0\ (\mathfrak{P});
\end{align*}
and so, by the properties of prime ideals,
\begin{align*}
\mathfrak{P} \equiv 0\ (\bar{\mathfrak{P}}); \quad \bar{\mathfrak{P}} \equiv 0\ (\mathfrak{P}); \quad \mathfrak{P} = \bar{\mathfrak{P}}.
\end{align*}
In conclusion, we have \\

\textbf{Theorem V.} \emph{For every primary ideal} $\mathfrak{Q}$ \emph{there exists one, and only one, prime ideal} $\mathfrak{P}$ \emph{which is a divisor of} $\mathfrak{Q}$ \emph{and a power of which is divisible by} $\mathfrak{Q}$; $\mathfrak{P}$ \emph{shall be referred to as the \textquotedblleft associated prime ideal"}.\footnote{The example $\mathfrak{M} = (x^{2}, xy)$ shows that the converse does not hold. The prime ideal $(x)$ satisfies all requirements, but $\mathfrak{M}$ is not primary.} $\mathfrak{P}$ \emph{is defined as the greatest common divisor of all ideals} $\mathfrak{B}$ \emph{with the property that a power of} $\mathfrak{B}$ \emph{is divisible by} $\mathfrak{Q}$. \emph{If} $\varrho$ \emph{is the smallest number such that} $\mathfrak{P}^{\varrho} \equiv 0\ (\mathfrak{Q})$, \emph{then} $\varrho$ \emph{shall be referred to as the exponent of} $\mathfrak{Q}$.\footnote{It is not in general true that $\mathfrak{P}^{\varrho} = \mathfrak{Q}$, as it is in the rings of the integers and the algebraic integers; for example:
\begin{align*}
\mathfrak{Q} = (x^{2}, y); \quad \mathfrak{P} = (x, y); \quad \mathfrak{P}^{2} = (x^{2}, xy, y^{2}) \equiv 0\ (\mathfrak{Q});\ \text{but}\ \mathfrak{Q} \not\equiv 0\ (\mathfrak{P}^{2});
\end{align*}
therefore $\mathfrak{Q}$ differs from $\mathfrak{P}^{2}$.} \\

We now prove, as the connection between primary and irreducible: \\

\textbf{Theorem VI.} \emph{Every non-primary ideal is reducible; in other words, every irreducible ideal is primary}.\footnote{The following example shows that the converse does not hold here: 
\begin{align*}
\mathfrak{Q} = (x^{2}, xy, y^{\lambda}) = [(x^{2}, y), (x, y^{\lambda})],
\end{align*}
where $\lambda \geq 2$. Here $\mathfrak{Q}$ is primary, but reducible. ($\mathfrak{Q}$ is primary because it contains all products of powers of $x$ and $y$ with a total dimension of $\lambda$; for every polynomial without a constant term, one power is therefore divisible by $\mathfrak{Q}$. However, should the polynomial $b$ in $a \cdot b \equiv 0\ (\mathfrak{Q})$ contain a constant term \--- so $b^{x} \not\equiv 0\ (\mathfrak{Q})$ for every $x$ \--- then $a$ must be divisible by $\mathfrak{Q}$, because every homogeneous component of $a \cdot b$ is divisible by $\mathfrak{Q}$ due to the homogeneity of the basis polynomials of $\mathfrak{Q}$.)} \\

Let $\mathfrak{K}$ be a \emph{non-primary} ideal, so that by Definition III there exists at least one pair of elements $a$, $b$ such that
\begin{align}
a \cdot b \equiv 0\ (\mathfrak{K}); \quad a \not\equiv 0\ (\mathfrak{K}); \quad b^{x} \not\equiv 0\ (\mathfrak{K}) \quad \text{for every } x.
\end{align}
We now construct the two ideals
\begin{align*}
\mathfrak{L}_{0} = (\mathfrak{K}, a); \quad \mathfrak{N}_{0} = (\mathfrak{K}, b),
\end{align*}
which by (3) are proper divisors of $\mathfrak{K}$, and by (3) it holds that
\begin{align}
\mathfrak{L}_{0} \cdot \mathfrak{N}_{0} \equiv 0\ (\mathfrak{K}).
\end{align}
For the elements $f$ of the least common multiple $\mathfrak{K}_{0} = [\mathfrak{L}_{0}, \mathfrak{N}_{0}]$, the following options now arise:

\emph{Either} from
\begin{align*}
f \equiv 0\ (\mathfrak{L}_{0}); \quad f \equiv 0\ (\mathfrak{N}_{0}),\ \text{that is } f \equiv a_{1} \cdot b\ (\mathfrak{K})
\end{align*}
always follows a representation
\begin{align*}
f \equiv l_{0} \cdot b\ (\mathfrak{K}); \quad l_{0} \equiv 0\ (\mathfrak{L}_{0});
\end{align*}
therefore, by (4), $f \equiv 0\ (\mathfrak{K})$, which implies $\mathfrak{K}_{0} \equiv 0\ (\mathfrak{K})$, and because $\mathfrak{K} \equiv 0\ (\mathfrak{K}_{0})$, it also holds that $\mathfrak{K} = \mathfrak{K}_{0}$, by which $\mathfrak{K}$ is proven to be \emph{reducible}.

\emph{Or} there is at least one $f \equiv 0\ (\mathfrak{K}_{0})$ for which there exists no such $l_{0}$. Using the $a_{1}$ belonging to this $f$ we then construct:
\begin{align*}
\mathfrak{L}_{1} = (\mathfrak{L}_{0}, a_{1}) = (\mathfrak{K}, a, a_{1}); \quad \mathfrak{N}_{1} = (\mathfrak{K}, b^{2}).
\end{align*}
Then, because $a_{1} \cdot b \equiv 0\ (\mathfrak{L}_{0})$ by (4), it also holds that
\begin{align*}
\mathfrak{L}_{1} \cdot \mathfrak{N}_{1} \equiv 0\ (\mathfrak{K}); \tag{4'}
\end{align*}
and $\mathfrak{L}_{1}$ is a \emph{proper} divisor of $\mathfrak{L}_{0}$.

For the elements $f$ of $\mathfrak{K}_{1} = [\mathfrak{L}_{1}, \mathfrak{N}_{1}]$, the same options occur:

\emph{Either} from
\begin{align*}
f \equiv 0\ (\mathfrak{L}_{1}); \quad f \equiv 0\ (\mathfrak{N}_{1}),\ \text{that is } f \equiv a_{2} \cdot b^{2}\ (\mathfrak{K})
\end{align*}
it always follows that
\begin{align*}
f \equiv l_{1} \cdot b^{2}; \quad l_{1} \equiv 0\ (\mathfrak{L}_{1});
\end{align*}
and so by (4'):
\begin{align*}
\mathfrak{K} = \mathfrak{K}_{1}.
\end{align*}

\emph{Or} there is at least one $f$ for which there exists no such $l_{1}$, which leads to the construction of $\mathfrak{L}_{2} = (\mathfrak{L}_{1}, a_{2})$, $\mathfrak{N}_{2} = (\mathfrak{K}, b^{2^{2}})$, with $\mathfrak{L}_{2} \cdot \mathfrak{N}_{2} \equiv 0\ (\mathfrak{K})$, where $\mathfrak{L}_{2}$ is a \emph{proper} divisor of $\mathfrak{L}_{1}$. Therefore, continuing in this way, in general we define
\begin{align*}
\mathfrak{L}_{0} &= (\mathfrak{K}, a);\ \mathfrak{L}_{1} = (\mathfrak{L}_{0}, a_{1});\ \dots;\ \mathfrak{L}_{\nu} = (\mathfrak{L}_{\nu - 1}, a_{\nu});\ \dots; \\
\mathfrak{N}_{0} &= (\mathfrak{K}, b);\ \mathfrak{N}_{1} = (\mathfrak{K}, b^{2});\ \dots;\ \mathfrak{N}_{\nu} = (\mathfrak{K}, b^{2^{\nu}});\ \dots,
\end{align*}
where the $a_{i}$ are defined so that there exists an $f$ such that
\begin{align*}
f \equiv 0\ (\mathfrak{L}_{i-1}); \quad f \equiv 0\ (\mathfrak{N}_{i-1}),
\end{align*}
that is
\begin{align*}
f \equiv a_{i} \cdot b^{2^{i-1}}\ (\mathfrak{K}); \quad \text{but } a_{i} \not\equiv 0\ (\mathfrak{L}_{i-1}).
\end{align*}
Subsequently it holds in general that $\mathfrak{L}_{i} \cdot \mathfrak{N}_{i} \equiv 0\ (\mathfrak{K})$, that by (3) $\mathfrak{N}_{i}$ is a \emph{proper} divisor of $\mathfrak{K}$, and that $\mathfrak{L}_{i}$ is a \emph{proper} divisor of $\mathfrak{L}_{i-1}$. By Theorem I of the Finite Chain, \emph{the chain of the} $\mathfrak{L}$ \emph{must therefore terminate in finitely many steps, say at} $\mathfrak{L}_{n}$. For each $f \equiv 0\ (\mathfrak{L}_{n})$, $f \equiv 0\ (\mathfrak{N}_{n})$, it follows that $f \equiv l_{n} \cdot b^{2^{n}}\ (\mathfrak{K})$, with $l_{n} \equiv 0\ (\mathfrak{L}_{n})$, and consequently by the above conclusion $\mathfrak{K} = [\mathfrak{L}_{n}, \mathfrak{N}_{n}]$, \emph{which proves that} $\mathfrak{K}$ \emph{is reducible}.

As a result of the preceding proofs, the \emph{uniqueness of the associated prime ideals} emerges as follows:

Let
\begin{align*}
\mathfrak{M} = [\mathfrak{B}_{1}, \dots, \mathfrak{B}_{k}] = [\mathfrak{D}_{1}, \dots, \mathfrak{D}_{k}]
\end{align*}
be two shortest, and therefore reduced, representations of $\mathfrak{M}$ as the least common multiple of irreducible ideals, the numbers of components of which are equal by Theorem IV. Then, by this theorem, the intermediate representations
\begin{align*}
\mathfrak{M} &= [\mathfrak{D}_{1}, \dots, \mathfrak{D}_{i-1}, \mathfrak{B}_{i}, \mathfrak{B}_{i+1}, \dots, \mathfrak{B}_{k}] = [\mathfrak{D}_{1}, \dots, \mathfrak{D}_{i-1}, \mathfrak{D}_{i}, \mathfrak{B}_{i+1}, \dots, \mathfrak{B}_{k}]
\\ &= [\bar{\mathfrak{A}}_{i}, \mathfrak{B}_{i}] = [\bar{\mathfrak{A}}_{i}, \mathfrak{D}_{i}]
\end{align*}
occurring there (where, as was remarked there, the index $j_{i} = i$ can be set) are also shortest representations. It therefore comes about that
\begin{align*}
\bar{\mathfrak{A}}_{i} \cdot \mathfrak{B}_{i} \equiv 0\ (\mathfrak{D}_{i}),\ \bar{\mathfrak{A}}_{i} \not\equiv 0\ (\mathfrak{D}_{i}); \quad \bar{\mathfrak{A}}_{i} \cdot \mathfrak{D}_{i} \equiv 0\ (\mathfrak{B}_{i}),\ \bar{\mathfrak{A}}_{i} \not\equiv 0\ (\mathfrak{B}_{i}).
\end{align*}
Because now by Theorem IV the irreducible ideals $\mathfrak{B}_{i}$ and $\mathfrak{D}_{i}$ are primary, it follows that there exist two numbers $\lambda_{i}$ and $\mu_{i}$ such that
\begin{align}
\mathfrak{B}_{i}^{\lambda_{i}} \equiv 0\ (\mathfrak{D}_{i}); \quad \mathfrak{D}_{i}^{\mu_{i}} \equiv 0\ (\mathfrak{B}_{i}).
\end{align}
Now let $\mathfrak{P}_{i}$ and $\bar{\mathfrak{P}}_{i}$ denote the corresponding prime ideals of $\mathfrak{B}_{i}$ and $\mathfrak{D}_{i}$ respectively; therefore $\mathfrak{P}_{i}^{\varrho_{i}} \equiv 0\ (\mathfrak{B}_{i})$ and $\bar{\mathfrak{P}}_{i}^{\sigma_{i}} \equiv 0\ (\mathfrak{D}_{i})$, and so by (5)
\begin{align*}
\mathfrak{P}_{i}^{\lambda_{i} \varrho_{i}} \equiv 0\ (\bar{\mathfrak{P}}_{i}); \quad \bar{\mathfrak{P}}_{i}^{\mu_{i} \sigma_{i}} \equiv 0\ (\mathfrak{P}_{i}),
\end{align*}
and from this, by the property of prime ideals:
\begin{align*}
\mathfrak{P}_{i} \equiv 0\ (\bar{\mathfrak{P}}_{i}); \quad \bar{\mathfrak{P}}_{i} \equiv 0\ (\mathfrak{P}_{i}); \quad \mathfrak{P}_{i} = \bar{\mathfrak{P}}_{i}.
\end{align*}
With this we have proven \\

\textbf{Theorem VII}. \emph{For two different shortest representations of an ideal as the least common multiple of irreducible ideals, the associated prime ideals agree, and in fact the same associated prime ideals}\footnote{This is shown, for instance, in the example from footnote 19:
\begin{align*}
(x^{2}, xy, y^{\lambda}) = [(x^{2}, y), (x, y^{\lambda})], 
\end{align*}
where $\lambda \geq 2$. Both corresponding prime ideals are $(x, y)$ here.} \emph{also occur for every decomposition thereof.} \emph{The ideals themselves can be paired up in at least one way such that a power of the ideal} $\mathfrak{B}_{i}$ \emph{is divisible by the associated} $\mathfrak{D}_{i}$ \emph{and vice versa. The numbers of ideals agree by Theorem} IV.\footnote{For the uniqueness of the \textquotedblleft isolated" ideals found among the irreducible ideals, see \S 7.} 

\section*{\S 5. Representation of an ideal as the least common multiple of maximal primary ideals. Uniqueness of the associated prime ideals.}

\tab \textbf{Definition IV}. \emph{A shortest representation} $\mathfrak{M} = [\mathfrak{Q}_{1}, \dots, \mathfrak{Q}_{\alpha}]$ \emph{is called the least common multiple of maximal primary ideals if all} $\mathfrak{Q}$ \emph{are primary, but the least common multiple of two} $\mathfrak{Q}$ \emph{is no longer primary}. \\

That at least one such representation always exists follows from the representation of $\mathfrak{M}$ as the least common multiple of irreducible ideals. This is because these ideals are primary; either there now already exists a representation through maximal primary ideals, or else the least common multiple of some two ideals is again primary. Because taking this least common multiple decreases the number of ideals by one, continuation of this procedure leads to the desired representation in finitely many steps.

This representation is reduced by Lemma IV. Conversely, every reduced representation arises through maximal primary ideals in this way, as the decomposition of the $\mathfrak{Q}$ into irreducible ideals shows.

In order to achieve an appropriate theorem of uniqueness here from Theorem VII, the connection with the associated prime ideals must be investigated, via \\

\textbf{Theorem VIII}. \emph{Should the primary ideals} $\mathfrak{N}_{1}, \mathfrak{N}_{2}, \dots, \mathfrak{N}_{\lambda}$ \emph{all have the same associated prime ideal} $\mathfrak{P}$, \emph{then their least common multiple} $\mathfrak{Q} = [\mathfrak{N}_{1}, \mathfrak{N}_{2}, \dots, \mathfrak{N}_{\lambda}]$ \emph{is also primary and has} $\mathfrak{P}$ \emph{as its associated prime ideal. If conversely} $\mathfrak{Q} = [\mathfrak{N}_{1}, \dots, \mathfrak{N}_{\lambda}]$ \emph{is a reduced representation for the primary ideal} $\mathfrak{Q}$, \emph{then all} $\mathfrak{N}_{i}$ \emph{are primary and have as their associated prime ideal the associated prime ideal} $\mathfrak{P}$ \emph{of} $\mathfrak{Q}$. \\

To prove the first part of the statement, first note that $\mathfrak{P}^{\varrho_{i}} \equiv 0\ (\mathfrak{N}_{i})$ for each $i$ also implies $\mathfrak{P}^{\tau} \equiv 0\ (\mathfrak{Q})$, where $\tau$ denotes the largest of the indices $\varrho_{i}$. Because $\mathfrak{P}$ is furthermore also a divisor of $\mathfrak{Q}$, $\mathfrak{P}$ is necessarily the associated prime ideal, if $\mathfrak{Q}$ is primary. It follows from
\begin{align*}
\mathfrak{A} \cdot \mathfrak{B} \equiv 0\ (\mathfrak{Q}); \quad \mathfrak{B}^{k} \not\equiv 0\ (\mathfrak{Q})\ \text{(for every $k$)}
\end{align*}
that consequently
\begin{align*}
\mathfrak{B} \not\equiv 0\ (\mathfrak{P});\ \text{so}\ \mathfrak{B}^{k} \not\equiv (\mathfrak{N}_{i})\ \text{(for every $k$)}
\end{align*}
and as a result
\begin{align*}
\mathfrak{A} \equiv 0\ (\mathfrak{N}_{i});\ \text{so}\ \mathfrak{A} \equiv 0\ (\mathfrak{Q}),
\end{align*}
which proves that $\mathfrak{Q}$ is primary, and that $\mathfrak{P}$ is the associated prime ideal.

Conversely, first let
\begin{align*}
\mathfrak{Q} = [\mathfrak{N}_{1}, \dots, \mathfrak{N}_{\lambda}] = [\mathfrak{N}_{i}, \mathfrak{L}_{i}]
\end{align*}
be a shortest representation of $\mathfrak{Q}$ (where $\mathfrak{Q}$ is primary) using \emph{primary} ideals $\mathfrak{N}_{i}$, and let $\mathfrak{P}_{i}$ be the respective associated prime ideals. It follows from
\begin{align*}
\mathfrak{L}_{i} \cdot \mathfrak{N}_{i} \equiv 0\ (\mathfrak{Q}); \quad \mathfrak{L}_{i} \not\equiv 0\ (\mathfrak{Q})
\end{align*}
(because of the shortest representation) that
\begin{align*}
\mathfrak{N}_{i}^{\sigma_{i}} \equiv 0\ (\mathfrak{Q});\ \text{or}\ \mathfrak{P}_{i}^{\varrho_{i} \sigma_{i}} \equiv 0\ (\mathfrak{Q}).
\end{align*}
Because the $\mathfrak{P}_{i}$ are all divisors of $\mathfrak{Q}$, $\mathfrak{P}_{i}$ is therefore equal to the associated prime ideal $\mathfrak{P}$ of $\mathfrak{Q}$ for all $i$.

It remains to show that for \emph{every reduced} representation $\mathfrak{Q} = [\mathfrak{N}_{1}, \dots, \mathfrak{N}_{\lambda}]$, the $\mathfrak{N}_{i}$ are primary.\footnote{The example $\mathfrak{Q} = [x^{2}, xy, y^{\lambda}] = [(x^{2}, xy, y^{2}, yz), (x, y^{\lambda})]$, where $\lambda \geq 2$, which is not reduced, shows that the reduced representation is crucial here. Here $(x^{2}, xy, y^{2}, yz) = [(x^{2}, y), (x, y^{2}, z)]$ is not primary by the above proof; this is because the last representation is a shortest one through primary ideals, but the associated prime ideals $(x, y)$ and $(x, y, z)$ are different. ($\mathfrak{Q}$ is primary by footnote 19.)} For this, decompose the $\mathfrak{N}_{i}$ into their irreducible ideals $\mathfrak{B}$; in the resulting \emph{shortest} representation $\mathfrak{Q} = [\mathfrak{B}_{1}, \dots, \mathfrak{B}_{\mu}]$, each ideal is then primary and has $\mathfrak{P}$ as its associated prime ideal by the previous parts of the proof. The same then holds for each $\mathfrak{N}_{i}$ by the first part of the theorem, proved above, which \emph{completes the proof of Theorem} VIII. \\

\textbf{Remark}. From this it follows that \emph{a prime ideal is necessarily irreducible}. This is because the reduced representation $\mathfrak{P} = [\mathfrak{N}, \mathfrak{L}]$ gives $\mathfrak{N} \equiv 0\ (\mathfrak{P})$ by Theorem VIII, and so, because $\mathfrak{P} \equiv 0\ (\mathfrak{N})$, this also gives $\mathfrak{P} = \mathfrak{N}$ and similarly $\mathfrak{P} = \mathfrak{L}$. The irreducibility of $\mathfrak{P}$ also follows directly, because $\mathfrak{P} = [\mathfrak{N}, \mathfrak{L}]$ implies $\mathfrak{N} \cdot \mathfrak{L} \equiv 0\ (\mathfrak{P})$, $\mathfrak{N} \not\equiv (\mathfrak{P})$, $\mathfrak{L} \not\equiv 0\ (\mathfrak{P})$, in contradiction with the defining property of a prime ideal. \\

Now let
\begin{align*}
\mathfrak{M} = [\mathfrak{Q}_{1}, \dots, \mathfrak{Q}_{\alpha}] = [\bar{\mathfrak{Q}}_{1}, \dots, \bar{\mathfrak{Q}}_{\beta}]
\end{align*}
be two reduced representations of $\mathfrak{M}$ as the least common multiple of maximal primary ideals. By decomposing the $\mathfrak{Q}$ into their irreducible ideals $\mathfrak{B}$ and the $\bar{\mathfrak{Q}}$ respectively into the irreducible ideals $\mathfrak{D}$, two reduced representations of $\mathfrak{M}$ as the least common multiple of irreducible ideals are produced, in which by Theorem VII both the numbers of components and the associated prime ideals are the same. By Theorem VIII, all irreducible ideals $\mathfrak{B}$ for a fixed $\mathfrak{Q}_{i}$ have \emph{the same associated prime ideal} $\mathfrak{P}_{i}$, while the $\mathfrak{P}_{k}$ associated with $\mathfrak{Q}_{k}$ is necessarily different from this, because otherwise by Theorem VIII no representation using maximal primary ideals exists. The number $\alpha$ of $\mathfrak{Q}$ is therefore equal to the number of different associated prime ideals $\mathfrak{P}$ of the $\mathfrak{B}$; these different $\mathfrak{P}$ construct the associated prime ideals of the $\mathfrak{Q}$. The same applies to the $\bar{\mathfrak{Q}}$ with respect to their decomposition into the $\mathfrak{D}$. From Theorem VII it therefore follows that \emph{the number of} $\mathfrak{Q}$ \emph{and} $\bar{\mathfrak{Q}}$ \emph{are the same}, and that \emph{their corresponding prime ideals are the same}. Together these show that the summary given at the start of the section about the irreducible ideals among maximal primary ideals holds true; that is, these, and only these, all have the same associated prime ideal. Theorem VII further shows the property of irreducibility of the maximal primary ideals: they admit no reduced representation as the least common multiple of maximal primary ideals.

In summary the following has been proved: \\

\textbf{Theorem IX}. \emph{For two reduced representations of an ideal as the least common multiple of maximal primary ideals, the numbers of components and the associated prime ideals (which are all different from each other) are the same. In other words, each} $\mathfrak{Q}$ \emph{can be uniquely associated with a} $\bar{\mathfrak{Q}}$ \emph{so that a power of} $\mathfrak{Q}$ \emph{is divisible by} $\bar{\mathfrak{Q}}$, \emph{and vice versa}.\footnote{One example of different representations is that given in footnote 12 for Theorem II: $(x^{2}, xy) = [(x), (x^{2}, \mu x + y)]$ for arbitrary $\mu$. Because the associated prime ideals $\mathfrak{P}_{1} = (x)$, $\mathfrak{P}_{2} = (x, y)$ are different, it consists of maximal primary ideals. For the uniqueness of the \textquotedblleft isolated" maximal primary ideals, see \S 7.} \emph{The $\mathfrak{Q}$ and $\bar{\mathfrak{Q}}$ have the property of irreducibility with respect to the decomposition into maximal primary ideals}. \\

\textbf{Remark}. Note that Theorem IX remains for the most part true if instead of reduced representations, only \emph{shortest} representations are required. If for instance $\mathfrak{M} = [\mathfrak{Q}_{1}, \dots, \mathfrak{Q}_{i}^{*}, \dots, \mathfrak{Q}_{\alpha}]$ is then reduced with respect to $\mathfrak{Q}_{i}^{*}$, and $\mathfrak{Q}_{i}^{*}$ is a proper divisor of $\mathfrak{Q}_{i}$, then by Lemma IV $\mathfrak{Q}_{i} = [\mathfrak{Q}_{i}^{*}, (\mathfrak{L}_{i}, \mathfrak{Q}_{i})]$, where $\mathfrak{L}_{i}$ is the complement of $\mathfrak{Q}_{i}$; this representation is reduced with respect to $\mathfrak{Q}_{i}^{*}$. By Theorem VIII, where in this application $(\mathfrak{L}_{i}, \mathfrak{Q}_{i})$ is replaced with a proper divisor if necessary, $\mathfrak{Q}_{i}^{*}$ is therefore primary and has the same associated prime ideal $\mathfrak{P}_{i}$ as $\mathfrak{Q}_{i}$. Continuation of this procedure shows that every such representation can be assigned a reduced representation through maximal primary ideals such that the number of components and the associated prime ideals are the same. \emph{It therefore also holds for shortest representations that in two different representations the number of components and the associated prime ideals are the same}.

The thus \emph{uniquely} defined associated prime ideals that are different from each other shall be called \emph{\textquotedblleft the associated prime ideals of $\mathfrak{M}$"} for short.

\section*{\S 6. Unique representation of an ideal as the least common multiple of relatively prime irreducible ideals.}

\tab \textbf{Definition V.} \emph{An ideal} $\mathfrak{R}$ \emph{is called relatively prime to} $\mathfrak{S}$ \emph{if} $\mathfrak{T} \cdot \mathfrak{R} \equiv 0\ (\mathfrak{S})$ \emph{necessarily implies} $\mathfrak{T} \equiv 0\ (\mathfrak{S})$. \emph{If} $\mathfrak{R}$ \emph{is relatively prime to} $\mathfrak{S}$ \emph{and} $\mathfrak{S}$ \emph{is also relatively prime to} $\mathfrak{R}$, \emph{then} $\mathfrak{R}$ \emph{and} $\mathfrak{S}$ \emph{are called mutually relatively prime}.\footnote{The relation of being relatively prime is not symmetric. For example, $\mathfrak{R} = (x^{2}, y)$ is relatively prime to $\mathfrak{S} = (x)$, but $\mathfrak{S}$ is not relatively prime to $\mathfrak{R}$, because $\mathfrak{S}^{2} \equiv 0\ (\mathfrak{R})$, whereas $\mathfrak{T} = \mathfrak{S} \not\equiv 0\ (\mathfrak{R})$.} \emph{An ideal is called relatively prime irreducible if it cannot be expressed as the least common multiple of mutually relatively prime proper divisors.} \\

In particular, if instead of $\mathfrak{T}$ the greatest common divisor $\mathfrak{T}_{0}$ of all $\mathfrak{T}$ for which $\mathfrak{T} \cdot \mathfrak{R} \equiv 0\ (\mathfrak{S})$ is taken, then we also have $\mathfrak{T}_{0} \cdot \mathfrak{R} \equiv 0\ (\mathfrak{S})$ and $\mathfrak{S} \equiv 0\ (\mathfrak{T}_{0})$. Therefore $\mathfrak{T}_{0} = \mathfrak{S}$ if $\mathfrak{R}$ is relatively prime to $\mathfrak{S}$, and $\mathfrak{T}_{0}$ is a proper divisor of $\mathfrak{S}$ if $\mathfrak{R}$ is not relatively prime to $\mathfrak{S}$.\footnote{This so defined $\mathfrak{T}_{0}$ is the same as Lasker's \textquotedblleft \underline{residual module}" and Dedekind's \textquotedblleft quotient" of two modules in his expansion into modules rather than ideals. Lasker, loc. cit. p49, Dedekind (Zahlentheorie), p504.}

The proof of uniqueness is underpinned by \\

\textbf{Theorem X.} 1. \emph{If $\mathfrak{R}$ is relatively prime to the ideals $\mathfrak{S}_{1}, \dots, \mathfrak{S}_{\lambda}$, then $\mathfrak{R}$ is also relatively prime to their least common multiple $\mathfrak{S}$.}

2. \emph{If the ideals $\mathfrak{S}_{1}, \dots, \mathfrak{S}_{\lambda}$ are relatively prime to $\mathfrak{R}$, then their least common multiple $\mathfrak{S}$ is also relatively prime to $\mathfrak{R}$.}

3. \emph{If $\mathfrak{R}$ is relatively prime to $\mathfrak{S}$ and $\mathfrak{S} = [\mathfrak{S}_{1}, \dots, \mathfrak{S}_{\lambda}]$ is a reduced representation of $\mathfrak{S}$, then $\mathfrak{R}$ is also relatively prime to each $\mathfrak{S}_{i}$.}

4. \emph{If $\mathfrak{S}$ is relatively prime to $\mathfrak{R}$, then each divisor $\mathfrak{S}_{i}$ of $\mathfrak{S}$ is also relatively prime to $\mathfrak{R}$.} \\

1. Because $\mathfrak{T} \cdot \mathfrak{R} \equiv 0\ (\mathfrak{S})$ necessarily implies $\mathfrak{T} \cdot \mathfrak{R} \equiv 0\ (\mathfrak{S}_{i})$, our assumption implies $\mathfrak{T} \equiv 0\ (\mathfrak{S}_{i})$, and therefore $\mathfrak{T} \equiv 0\ (\mathfrak{S})$.

2. Let $\mathfrak{C}_{1} = [\mathfrak{S}_{2}, \dots, \mathfrak{S}_{\lambda}]$, $\mathfrak{C}_{12} = [\mathfrak{S}_{3}, \dots, \mathfrak{S}_{\lambda}]$, $\dots$, $\mathfrak{C}_{12 \dots \lambda - 1} = \mathfrak{S}_{\lambda}$. From $\mathfrak{T} \cdot \mathfrak{S} \equiv 0\ (\mathfrak{R})$ it then follows that $\mathfrak{T} \cdot \mathfrak{C}_{1} \cdot \mathfrak{S}_{1} \equiv 0\ (\mathfrak{R})$, and so by our assumption $\mathfrak{T} \cdot \mathfrak{C}_{1} \equiv 0\ (\mathfrak{R})$. From this it follows furthermore that $\mathfrak{T} \cdot \mathfrak{C}_{12} \cdot \mathfrak{S}_{2} \equiv 0\ (\mathfrak{R})$, and so by our assumption $\mathfrak{T} \cdot \mathfrak{C}_{12} \equiv 0\ (\mathfrak{R})$. Proceeding in this way, it follows finally that $\mathfrak{T} \cdot \mathfrak{C}_{12 \dots \lambda - 1} = \mathfrak{T} \cdot \mathfrak{S}_{\lambda} \equiv 0\ (\mathfrak{R})$, and so $\mathfrak{T} \equiv 0\ (\mathfrak{R})$.

3. Let $\mathfrak{C}_{i}$ denote the complement of $\mathfrak{S}_{i}$; it then follows from $\mathfrak{T}_{0} \cdot \mathfrak{R} \equiv 0\ (\mathfrak{S}_{i})$ that $[\mathfrak{T}_{0}, \mathfrak{C}_{i}] \cdot \mathfrak{R} \equiv 0\ (\mathfrak{S})$, since $\mathfrak{C}_{i} \cdot \mathfrak{R} \equiv 0\ (\mathfrak{C}_{i})$.

Should $\mathfrak{T}_{0}$ now be a proper divisor of $\mathfrak{S}_{i}$, then because of the reduced representation $\mathfrak{S} = [\mathfrak{S}_{i}, \mathfrak{C}_{i}]$, it is also true that $[\mathfrak{T}_{0}, \mathfrak{C}_{i}]$ is a proper divisor of $\mathfrak{S}$, contradicting our assumption.

4. It follows from $\mathfrak{T}_{0} \cdot \mathfrak{S}_{i} \equiv 0\ (\mathfrak{R})$, where $\mathfrak{T}_{0}$ is a proper divisor of $\mathfrak{R}$, that $\mathfrak{T}_{0} \cdot \mathfrak{S} \equiv 0\ (\mathfrak{R})$; this contradicts our assumption. \\ 

Definition V shows in particular that every ideal is relatively prime to the \emph{trivial ideal} $\mathfrak{O}$ consisting of \emph{all} elements of $\Sigma$;\footnote{$\mathfrak{O}$ plays the role of the trivial ideal only with respect to divisibility and least common multiples, not with respect to the formation of products. For example, $\mathfrak{O} = (x)$ for the ring of all polynomials in $x$ with integer coefficients and without a constant term; $\mathfrak{O} = (2)$ for the ring of all even numbers.} the following theorems however apply only for ideals other than $\mathfrak{O}$.

From Theorem X, which has just been proved, it immediately follows that: \\

\textbf{Lemma V}. \emph{Every representation of an ideal as the least common multiple of mutually relatively prime ideals other than $\mathfrak{O}$ is reduced.} \\

Let $\mathfrak{M} = [\mathfrak{R}_{1}, \mathfrak{R}_{2}, \dots, \mathfrak{R}_{\sigma}] = [\mathfrak{R}_{i}, \mathfrak{L}_{i}]$ be one such representation. By Theorem X, part 2, $\mathfrak{L}_{i}$ is then also relatively prime to $\mathfrak{R}_{i}$; $\mathfrak{R}_{i}$ can therefore not appear in $\mathfrak{L}_{i}$, and so the representation is a shortest one. This is because if $\mathfrak{L}_{i} \equiv 0\ (\mathfrak{R}_{i})$, where $\mathfrak{L}_{i}$ is relatively prime to $\mathfrak{R}_{i}$, it would also follow that $\mathfrak{O} \equiv 0\ (\mathfrak{R}_{i})$, as $\mathfrak{O} \cdot \mathfrak{L}_{i} \equiv 0\ (\mathfrak{R}_{i})$, and it therefore holds that $\mathfrak{R}_{i} = \mathfrak{O}$, which is by assumption impossible. Now suppose $\mathfrak{R}_{i}$ can be replaced with the proper divisor $\mathfrak{R}_{i}^{*}$. Then $\mathfrak{R}_{i}$ is reducible by Lemma II, and $\mathfrak{R}_{i} = [\mathfrak{R}_{i}^{*}, (\mathfrak{R}_{i}, \mathfrak{L}_{i})]$. Replace $(\mathfrak{R}_{i}, \mathfrak{L}_{i})$ here with a proper divisor $(\mathfrak{R}_{i}, \mathfrak{L}_{i})^{*}$ if necessary, and likewise replace $\mathfrak{R}_{i}^{*}$ with a proper divisor if necessary so that a reduced representation for $\mathfrak{R}_{i}$ is formed. Then Theorem X, part 3, shows that $\mathfrak{L}_{i}$ is also relatively prime to $(\mathfrak{R}_{i}, \mathfrak{L}_{i})^{*}$, and this is different from $\mathfrak{O}$ because of the shortest representation. Because however $(\mathfrak{R}_{i}, \mathfrak{L}_{i})^{*}$ appears in $(\mathfrak{R}_{i}, \mathfrak{L}_{i})$ and $(\mathfrak{R}_{i}, \mathfrak{L}_{i})$ appears in $\mathfrak{L}_{i}$, a contradiction therefore arises. \\

From Theorem X a further theorem arises, concerning the connection with the associated prime ideals: \\

\textbf{Theorem XI}. \emph{If $\mathfrak{R}$ is relatively prime to $\mathfrak{S}$ and $\mathfrak{S}$ is different from $\mathfrak{O}$, then no associated prime ideal of $\mathfrak{R}$\footnote{The definition of associated prime ideals of $\mathfrak{R}$ is given in the conclusion of \S 5.} is divisible by an associated prime ideal of $\mathfrak{S}$. Conversely, should no such divisibility occur, then $\mathfrak{R}$ is relatively prime to $\mathfrak{S}$, and of course $\mathfrak{S}$ is different from $\mathfrak{O}$.} \\

Let
\begin{align*}
\mathfrak{R} = [\mathfrak{Q}_{1}, \dots, \mathfrak{Q}_{\alpha}], \quad \mathfrak{S} = [\mathfrak{Q}_{1}^{*}, \dots, \mathfrak{Q}_{\beta}^{*}]
\end{align*}
be reduced representations of $\mathfrak{R}$ and $\mathfrak{S}$ through maximal primary ideals, and let $\mathfrak{P}_{1}, \dots, \mathfrak{P}_{\alpha}, \mathfrak{P}_{1}^{*}, \dots, \mathfrak{P}_{\beta}^{*}$ be the associated prime ideals. We prove the statement in the following form: if a $\mathfrak{P}$ is divisible by a $\mathfrak{P}^{*}$, then $\mathfrak{R}$ cannot be relatively prime to $\mathfrak{S}$, and vice versa.

Therefore let $\mathfrak{P}_{\mu} \equiv 0\ (\mathfrak{P}_{\nu}^{*})$ and as a result also $\mathfrak{Q}_{\mu} \equiv 0\ (\mathfrak{P}_{\nu}^{*})$, from which by definition of $\mathfrak{P}_{\nu}^{*}$ it follows that $\mathfrak{Q}_{\mu}^{\sigma_{\nu}} \equiv 0\ (\mathfrak{Q}_{\nu}^{*})$, and therefore also $\mathfrak{R}^{\sigma_{\nu}} \equiv 0\ (\mathfrak{Q}_{\nu}^{*})$. Now let $\mathfrak{R}^{\tau}$ be the lowest power of $\mathfrak{R}$ divisible by $\mathfrak{Q}_{\nu}^{*}$. If $\tau = 1$, then $\mathfrak{R}$ is divisible by $\mathfrak{Q}_{\nu}^{*}$; because, however, $\mathfrak{S}$ is different from $\mathfrak{O}$ by assumption, $\mathfrak{R}$ is therefore not relatively prime to $\mathfrak{Q}_{\nu}^{*}$. If $\tau \geq 2$, then $\mathfrak{R}^{\tau - 1} \cdot \mathfrak{R} \equiv 0\ (\mathfrak{Q}_{\nu}^{*})$, $\mathfrak{R}^{\tau - 1} \not\equiv 0\ (\mathfrak{Q}_{\nu}^{*})$, and so $\mathfrak{R}$ is not relatively prime to $\mathfrak{Q}_{\nu}^{*}$,\footnote{$\mathfrak{R}^{0}$ is not defined, because $\Sigma$ does not need to contain a unit; therefore the case $\tau = 1$ must be considered separately. $\tau = 0$ is also impossible when $\Sigma$ does contain a unit, due to the assumption regarding $\mathfrak{S}$.} and hence in both cases $\mathfrak{R}$ is also not relatively prime to $\mathfrak{S}$ by Theorem X, part 3.

If conversely $\mathfrak{R}$ is not relatively prime to $\mathfrak{S}$, then by Theorem X, part 1, $\mathfrak{R}$ is also not relatively prime to at least one $\mathfrak{Q}_{\nu}^{*}$. It therefore holds that $\mathfrak{T}_{0} \cdot \mathfrak{R} \equiv 0\ (\mathfrak{Q}_{\nu}^{*})$, $\mathfrak{T}_{0} \not\equiv 0\ (\mathfrak{Q}_{\nu}^{*})$, and from this, because $\mathfrak{Q}_{\nu}^{*}$ is primary, it follows that $\mathfrak{R}^{\tau} \equiv 0\ (\mathfrak{Q}_{\nu}^{*})$, and so also $\mathfrak{Q}_{1}^{\tau} \dots \mathfrak{Q}_{\alpha}^{\tau} \equiv 0\ (\mathfrak{Q}_{\nu}^{*})$. It then follows, however, that $\mathfrak{P}_{1}^{\tau \varrho_{1}} \dots \mathfrak{P}_{\alpha}^{\tau \varrho_{\alpha}} \equiv 0\ (\mathfrak{P}_{\nu}^{*})$ for the associated prime ideals, and by the properties of the prime ideals, $\mathfrak{P}_{\nu}^{*}$ is therefore contained in at least one $\mathfrak{P}$, which \emph{completes the proof of Theorem} XI. \\

The \emph{existence\footnote{The existence of the decomposition can also be proved directly, in exact analogy with the proof of the existence of the decomposition into finitely many irreducible ideals given in \S 2.} and uniqueness of the decomposition into relatively prime irreducible ideals} now arises from Theorems X and XI as follows:

Let $\mathfrak{M} = [\mathfrak{Q}_{1}, \dots, \mathfrak{Q}_{\alpha}]$ be a reduced (or at least shortest) representation of $\mathfrak{M}$ through maximal primary ideals, and let $\mathfrak{P}_{1}, \dots, \mathfrak{P}_{\alpha}$ be the associated prime ideals. We collect the $\mathfrak{P}$ together in groups such that \emph{no ideal in a group is divisible by an ideal in a different group, and each individual group cannot be split into two subgroups that both have this property}. In order to construct such a grouping, it must be noted that by definition the group $G$ containing each ideal $\mathfrak{P}$ must also contain all its divisors and multiples (that is to say, divisible by $\mathfrak{P}$) occurring in $\mathfrak{P}_{1} \dots \mathfrak{P}_{\alpha}$. For example, let $\mathfrak{P}^{(i_{1})}$ be all the multiples of $\mathfrak{P}$, $\mathfrak{P}_{j_{1}}^{(i_{1})}$ all the divisors of $\mathfrak{P}^{(i_{1})}$, $\mathfrak{P}_{j_{1}}^{(i_{1}, i_{2})}$ all the multiples of $\mathfrak{P}_{j_{1}}^{(i_{1})}$ and so on; so in general $\mathfrak{P}_{j_{1} \dots j_{\lambda - 1}}^{(i_{1} \dots i_{\lambda})}$ is all the multiples of $\mathfrak{P}_{j_{1} \dots j_{\lambda - 1}}^{(i_{1} \dots i_{\lambda - 1})}$ and $\mathfrak{P}_{j_{1} \dots j_{\lambda}}^{(i_{1} \dots i_{\lambda})}$ is all the divisors of $\mathfrak{P}_{j_{1} \dots j_{\lambda - 1}}^{(i_{1} \dots i_{\lambda})}$. Because it deals with only finitely many ideals $\mathfrak{P}$, this algorithm must terminate in finitely many steps; that is, no ideal different from all preceding ones is left over as a result. The thus obtained system of ideals $\mathfrak{P}$ now in fact forms a group $G$ with the desired properties. 

By definition $G$ contains in addition to each $\mathfrak{P}_{j_{1} \dots j_{\lambda - 1}}^{(i_{1} \dots i_{\lambda - 1})}$ all multiples $\mathfrak{P}_{j_{1} \dots j_{\lambda - 1}}^{(i_{1} \dots i_{\lambda})}$, and in addition to each $\mathfrak{P}_{j_{1} \dots j_{\lambda - 1}}^{(i_{1} \dots i_{\lambda})}$ all divisors $\mathfrak{P}_{j_{1} \dots j_{\lambda}}^{(i_{1} \dots i_{\lambda})}$. Among these however are also contained all divisors of the $\mathfrak{P}_{j_{1} \dots j_{\lambda - 1}}^{(i_{1} \dots i_{\lambda - 1})}$, while the multiples of the $\mathfrak{P}_{j_{1} \dots j_{\lambda - 1}}^{(i_{1} \dots i_{\lambda})}$ are themselves again $\mathfrak{P}_{j_{1} \dots j_{\lambda - 1}}^{(i_{1} \dots i_{\lambda})}$. The ideals not contained in $G$ can therefore neither be divisors nor multiples of those ideals contained in $G$. $G$ however also satisfies the irreducibility condition. Because if a division into two subgroups $G^{(1)}$ and $G^{(2)}$ exists, and $G^{(1)}$ contains for instance $\mathfrak{P}_{j_{1} \dots j_{\lambda}}^{(i_{1} \dots i_{\lambda})}$ (respectively $\mathfrak{P}_{j_{1} \dots j_{\lambda - 1}}^{(i_{1} \dots i_{\lambda})}$), it then also contains all preceding ideals, because these are alternating multiples and divisors (respectively divisors and multiples). So it also contains $\mathfrak{P}$, and thus the whole group $G$. Should we continue accordingly with the ideals not contained in $G$, then a division of all $\mathfrak{P}$ into groups $G_{1}, \dots, G_{\sigma}$ which possess the desired properties is obtained. \emph{Such a grouping is unique}; because if $G_{1}^{'}, \dots, G_{\tau}^{'}$ is a second grouping and $\mathfrak{P}_{j_{1} \dots j_{\lambda}}^{(i_{1} \dots i_{\lambda})}$ (respectively $\mathfrak{P}_{j_{1} \dots j_{\lambda - 1}}^{(i_{1} \dots i_{\lambda})}$) is an element of $G_{i}^{'}$, then by the above reasoning $G_{i}^{'}$ contains the whole group $G$ and is therefore identically equal to $G$ because of the irreducibility property.

Now denote by $\mathfrak{P}_{i \mu}$ (where $\mu$ runs from $1$ to $\lambda_{i}$) the ideals $\mathfrak{P}$ collected in a group $G_{i}$. Then, because the $\mathfrak{P}$ are all different from each other, the associated primary ideals $\mathfrak{Q}_{i \mu}$ of a shortest representation are uniquely determined. Should we set
\begin{align*}
\mathfrak{R}_{i} = [\mathfrak{Q}_{i1}, \dots, \mathfrak{Q}_{i \lambda_{i}}],\ \text{then}\ \mathfrak{M} = [\mathfrak{R}_{1}, \dots, \mathfrak{R}_{\sigma}];
\end{align*}
we show that, as a result, a \emph{decomposition of $\mathfrak{M}$ into relatively prime irreducible ideals} is achieved. First of all, it must be noted that Theorem XI also remains applicable when $\mathfrak{R}_{i} = [\mathfrak{Q}_{i1}, \dots, \mathfrak{Q}_{i \lambda_{i}}]$ is only a shortest representation, because the associated prime ideals are uniquely defined, as seen from the remark on Theorem IX. Because now no associated prime ideal $\mathfrak{P}_{i \mu}$ of $\mathfrak{R}_{i}$  is divisible by an associated prime ideal $\mathfrak{P}_{j \nu}$ of $\mathfrak{R}_{j}$ and vice versa, by Theorem XI $\mathfrak{R}_{i}$ and $\mathfrak{R}_{j}$ are mutually relatively prime, and each individual $\mathfrak{R}_{i}$ is relatively prime irreducible by Theorem XI due to the irreducibility property of the group $G_{i}$, because with reducibility, ideals completely different from $\mathfrak{O}$ come into question. Furthermore, by Lemma V it is a reduced representation also if we had originally only started out with a shortest representation through the $\mathfrak{Q}$. Conversely, every decomposition into relatively prime irreducible ideals leads to the given grouping of the $\mathfrak{P}_{i \mu}$ by Theorem XI.

Now let $\mathfrak{M} = [\bar{\mathfrak{R}}_{1}, \dots, \bar{\mathfrak{R}}_{\tau}]$ be a \emph{second representation of $\mathfrak{M}$ through relatively prime irreducible ideals}, which is reduced by Lemma V. Because then, as the decomposition of the $\bar{\mathfrak{R}}$ into maximal primary ideals shows, the associated prime ideals are the same, the groupings of these prime ideals, the uniqueness of which was proven above, are therefore also the same. It therefore follows that $\tau = \sigma$; and the notation can be chosen so that $\mathfrak{R}_{i}$ and $\bar{\mathfrak{R}}_{i}$ belong to the same group. Therefore let
\begin{align*}
\mathfrak{M} = [\mathfrak{R}_{i}, \mathfrak{L}_{i}] = [\bar{\mathfrak{R}}_{i}, \bar{\mathfrak{L}}_{i}]
\end{align*}
be the representation through ideal and complement. Then, because $\bar{\mathfrak{R}}_{i}$ is associated with the same group as $\mathfrak{R}_{i}$, by Theorem XI $\mathfrak{L}_{i}$ is also \emph{relatively prime to} $\bar{\mathfrak{R}}_{i}$ and $\bar{\mathfrak{L}}_{i}$ \emph{relatively prime to} $\mathfrak{R}_{i}$. Because
\begin{align*}
\mathfrak{R}_{i} \cdot \mathfrak{L}_{i} \equiv 0\ (\bar{\mathfrak{R}}_{i}), \quad \bar{\mathfrak{R}}_{i} \cdot \bar{\mathfrak{L}}_{i} \equiv 0\ (\mathfrak{R}_{i})
\end{align*}
it therefore follows that
\begin{align*}
\mathfrak{R}_{i} \equiv 0\ (\bar{\mathfrak{R}}_{i}); \quad \bar{\mathfrak{R}}_{i} \equiv 0\ (\mathfrak{R}_{i}); \quad \mathfrak{R}_{i} = \bar{\mathfrak{R}}_{i}.
\end{align*}
Hence we have proved \\

\textbf{Theorem XII}. \emph{Every ideal can be uniquely expressed as the least common multiple of finitely many mutually relatively prime and relatively prime irreducible ideals.}

\section*{\S 7. Uniqueness of the isolated ideals.}

\tab \textbf{Definition VI}. \emph{If the shortest representation $\mathfrak{M} = [\mathfrak{R}, \mathfrak{L}]$ is reduced with respect to $\mathfrak{L}$, then $\mathfrak{R}$ is called an isolated ideal if no associated prime ideal of $\mathfrak{R}$ appears in an associated prime ideal of $\mathfrak{L}$, or in other words, if $\mathfrak{L}$ is relatively prime to $\mathfrak{R}$.} \\

Subsequently the representation $\mathfrak{M} = [\mathfrak{R}, \mathfrak{L}]$ satisfies the requirements of the representation $\mathfrak{M} = [\mathfrak{R}_{i}, \mathfrak{L}_{i}]$ in Lemma V, and with Lemma V we have proved \\

\textbf{Lemma VI}. \emph{If $\mathfrak{R}$ is an isolated ideal of the shortest and reduced (with respect to $\mathfrak{L}$) representation $\mathfrak{M} = [\mathfrak{R}, \mathfrak{L}]$, then the representation is also reduced with respect to $\mathfrak{R}$.} \\

Because we therefore always have a reduced representation when using isolated ideals, the associated prime ideals occurring in the decompositions of $\mathfrak{R}$ and $\mathfrak{L}$ into \emph{irreducible} ideals complement each other to give the uniquely determined associated prime ideals occurring in the corresponding decomposition of $\mathfrak{M}$. 

Therefore no prime ideal of $\mathfrak{R}$ belonging to the decomposition into irreducible ideals appears in the remaining prime ideals belonging to the corresponding decomposition of $\mathfrak{M}$. If conversely this condition is satisfied, and if $\mathfrak{R}$ features in at least one representation $\mathfrak{M} = [\mathfrak{R}, \mathfrak{L}]$ reduced with respect to $\mathfrak{R}$, and so also in a reduced representation $\mathfrak{M} = [\mathfrak{R}, \mathfrak{L}^{*}]$, then $\mathfrak{R}$ is isolated by Definition VI. From this the following definition arises, which is independent of the particular complement $\mathfrak{L}$: \\

\textbf{Definition VIa}. \emph{$\mathfrak{R}$ is called an isolated ideal if the prime ideals of $\mathfrak{R}$ belonging to the decomposition into irreducible ideals do not appear in the remaining prime ideals belonging to the corresponding decomposition of $\mathfrak{M}$, and if $\mathfrak{R}$ appears in at least one representation $\mathfrak{M} = [\mathfrak{R}, \mathfrak{L}]$ reduced with respect to $\mathfrak{R}$.\footnote{Should ordinary associated prime ideals (conclusion of \S 5) be introduced, then it would be added as a special requirement that those of $\mathfrak{L}$ all be different from those of $\mathfrak{R}$. Following Definition VIa the representation therefore needs no longer be assumed to be reduced with respect to the complement.}} \\

Now let
\begin{align*}
\mathfrak{M} = [\mathfrak{R}, \mathfrak{L}] = [\bar{\mathfrak{R}}, \bar{\mathfrak{L}}]
\end{align*}
be two representations of $\mathfrak{M}$ through isolated ideals $\mathfrak{R}$ and $\bar{\mathfrak{R}}$ and complements $\mathfrak{L}$ and $\bar{\mathfrak{L}}$ such that the associated prime ideals of $\mathfrak{R}$ and $\bar{\mathfrak{R}}$ are the same. Should $\mathfrak{L}$ and $\bar{\mathfrak{L}}$ be replaced with divisors $\mathfrak{L}^{*}$ and $\bar{\mathfrak{L}}^{*}$ such that the representations are reduced, then the associated prime ideals of $\mathfrak{L}^{*}$ and $\bar{\mathfrak{L}}^{*}$ are also the same; by Theorem XI $\mathfrak{L}^{*}$ is therefore relatively prime to $\bar{\mathfrak{R}}$ and $\bar{\mathfrak{L}}^{*}$ relatively prime to $\mathfrak{R}$. Because
\begin{align*}
\mathfrak{R} \cdot \mathfrak{L}^{*} \equiv 0\ (\bar{\mathfrak{R}}); \quad \bar{\mathfrak{R}} \cdot \bar{\mathfrak{L}}^{*} \equiv 0\ (\mathfrak{R})
\end{align*}
it therefore follows that
\begin{align*}
\mathfrak{R} \equiv 0\ (\bar{\mathfrak{R}}); \quad \bar{\mathfrak{R}} \equiv 0\ (\mathfrak{R}); \quad \mathfrak{R} = \bar{\mathfrak{R}};
\end{align*}
isolated ideals are therefore \emph{uniquely} determined by the associated prime ideals. This in particular results in a strengthening of Theorems VII and IX regarding the decomposition into irreducible and maximal primary ideals, where due to the remark on Theorem IX only shortest representations need to be assumed. In summary: \\

\textbf{Theorem XIII}. \emph{For each shortest representation of an ideal as the least common multiple of irreducible ideals (respectively maximal primary ideals), the isolated irreducible ideals (respectively maximal primary ideals) are uniquely determined; the non-uniqueness applies only to the non-isolated irreducible ideals (respectively maximal primary ideals).\footnote{This theorem is already given without proof by Macaulay for polynomial ideals in the case of the decomposition into maximal primary ideals; his definition of the isolated and non-isolated (imbedded) primary ideals can be viewed as the irrational version of that given below.} In general, the isolated ideals are uniquely determined by the associated prime ideals.} \\

If the ideals $\mathfrak{B}_{i}$ (respectively $\mathfrak{D}_{j}$) in one such shortest representation through irreducible ideals (respectively maximal primary ideals) are non-isolated, then by definition the complements $\mathfrak{A}_{i}$ (respectively $\mathfrak{L}_{j}$) are divisible by $\mathfrak{P}_{i}$ (respectively $\mathfrak{P}_{j}$). It therefore follows that
\begin{align*}
\mathfrak{A}_{i}^{\varrho_{i}} \equiv 0\ (\mathfrak{B}_{i}) \quad \text{(respectively } \mathfrak{L}_{j}^{\sigma_{j}} \equiv 0\ (\mathfrak{D}_{j}) \text{)}.
\end{align*}
It follows conversely from the fulfilment of these relations that $\mathfrak{P}_{i}$ (respectively $\mathfrak{P}_{j}$) appears in at least one associated prime ideal of the complement, and so by the remark on Theorem IX also in an associated prime ideal of the divisor $\mathfrak{L}_{j}^{*}$ of $\mathfrak{L}_{j}$ that leads to a representation that is reduced with respect to $\mathfrak{L}_{j}^{*}$; the $\mathfrak{B}_{i}$ (respectively $\mathfrak{D}_{j}$) are therefore non-isolated. In particular, irreducible ideals $\mathfrak{B}_{i}$ for which the associated prime ideal appears more than once in the decomposition of $\mathfrak{M}$ are always non-isolated. \emph{Non-isolated primary ideals are therefore also characterised by the fact that a power of each complement is divisible by them; isolated primary ideals are characterised by the fact that this cannot be satisfied.}

\section*{\S 8. Unique representation of an ideal as the product of coprime irreducible ideals.}

\tab Should the ring $\Sigma$ contain a unit, that is, an element $\varepsilon$ such that $\varepsilon \cdot a = a$ for all elements in $\Sigma$,\footnote{$\Sigma$ clearly cannot contain more than \emph{one} unit because of the commutativity of multiplication, since for some two units $\varepsilon_{1}$ and $\varepsilon_{2}$ it holds that $\varepsilon_{1} \varepsilon_{2} = \varepsilon_{2} = \varepsilon_{1}$.} then the coprime ideals can be defined by \\

\textbf{Definition VIII}. \emph{Two ideals $\mathfrak{R}$ and $\mathfrak{S}$ are called coprime if their greatest common divisor is the trivial ideal $\mathfrak{O} = (\varepsilon)$ consisting of all elements of $\Sigma$. An ideal is called coprime irreducible if it cannot be expressed as the least common multiple of pairwise coprime ideals.} \\ 

Note that two coprime ideals are always mutually relatively prime. By definition there are two elements $r \equiv 0\ (\mathfrak{R})$ and $s \equiv 0\ (\mathfrak{S})$ such that $\varepsilon = r + s$. It follows from $\mathfrak{T} \cdot \mathfrak{R} \equiv 0\ (\mathfrak{S})$ however that $\mathfrak{T} \cdot r \equiv 0\ (\mathfrak{S})$, and so $\mathfrak{T} \cdot \varepsilon = \mathfrak{T} \equiv 0\ (\mathfrak{S})$; similarly $\bar{\mathfrak{T}} \cdot \mathfrak{S} \equiv 0\ (\mathfrak{R})$ gives $\bar{\mathfrak{T}} \equiv 0\ (\mathfrak{R})$.\footnote{The converse does not hold, however; for example, the ideals $\mathfrak{R} = (x)$ and $\mathfrak{S} = (y)$ are mutually relatively prime, but not coprime.} It therefore follows from Lemma V that each representation through pairwise coprime ideals is \emph{reduced}.

The following theorem, analogous to Theorem X, lays the foundation for the proof of uniqueness: \\

\textbf{Theorem XIV}. \emph{If $\mathfrak{R}$ is coprime to each of the ideals $\mathfrak{S}_{1}, \dots, \mathfrak{S}_{\lambda}$, then $\mathfrak{R}$ is also coprime to $\mathfrak{S} = [\mathfrak{S}_{1}, \dots, \mathfrak{S}_{\lambda}]$. Conversely, it also follows from the coprimality of $\mathfrak{R}$ and $\mathfrak{S}$ that $\mathfrak{R}$ is coprime to each $\mathfrak{S}_{j}$. If $\mathfrak{R} = [\mathfrak{R}_{1}, \dots, \mathfrak{R}_{\mu}]$ and each $\mathfrak{R}_{i}$ is coprime to each $\mathfrak{S}_{j}$, then $\mathfrak{R}$ and $\mathfrak{S}$ are coprime; the converse again holds here.} \\

If $\mathfrak{R}$ is coprime to each $\mathfrak{S}_{j}$, then there exist elements $s_{j}$ such that
\begin{align*}
s_{j} \equiv 0\ (\mathfrak{S}_{j}); \quad s_{j} \equiv \varepsilon\ (\mathfrak{R}).
\end{align*}
Hence
\begin{align*}
s_{1} \cdot s_{2} \cdot \dots{} \cdot s_{\lambda} \equiv 0\ (\mathfrak{S}); \quad s_{1} \cdot s_{2} \cdot \dots{} \cdot s_{\lambda} \equiv \varepsilon\ (\mathfrak{R}), \quad (\mathfrak{R}, \mathfrak{S}) = (\varepsilon).
\end{align*}
Because, however, $(\mathfrak{R}, \mathfrak{S})$ is divisible by each $(\mathfrak{R}, \mathfrak{S}_{j})$, the converse also holds. Repeated application of this conclusion results in the second part of the statement. Namely, if $\mathfrak{R}_{i}$ is coprime to $\mathfrak{S}_{1}, \dots, \mathfrak{S}_{\lambda}$ for fixed $i$, then $\mathfrak{R}_{i}$ is coprime to $\mathfrak{S}$. If this holds for every $i$, then because the relationship of coprimality is symmetric, $\mathfrak{S}$ is coprime to $\mathfrak{R}$. Conversely, the coprimality of $\mathfrak{S}$ to $\mathfrak{R}_{i}$ follows from the coprimality of $\mathfrak{R}$ and $\mathfrak{S}$, and from this follows the coprimality of $\mathfrak{R}_{i}$ to $\mathfrak{S}_{j}$.

The proof of the existence and uniqueness\footnote{The existence of the decomposition can again be proved in direct analogy with \S 2; the proof of uniqueness can also be conducted directly (cf. that mentioned in the introduction about Schmeidler and Noether-Schmeidler). The proof given here also gives an insight into the structure of the coprime irreducible ideals.} of the decomposition into coprime irreducible ideals comes from a \emph{unique} grouping, like the corresponding proof for relatively prime ideals. Because, however, the relatively prime irreducible ideals $\mathfrak{R}_{1}, \dots, \mathfrak{R}_{\sigma}$ of $\mathfrak{M}$ are \emph{uniquely} defined by Theorem XII, referring back to the associated prime ideals is unnecessary here.

We collect together the uniquely defined relatively prime irreducible ideals $\mathfrak{R}_{1}, \dots, \mathfrak{R}_{\sigma}$ of $\mathfrak{M}$ in groups in such a way that \emph{each ideal of a group is coprime to every ideal of a group different from it, while each individual group cannot be split into two subgroups such that each ideal of a subgroup is coprime to every ideal of the other subgroup}. One such grouping is given as follows: by definition, for every ideal $\mathfrak{R}$ contained in each individual group $G$, all ideals not coprime to $\mathfrak{R}$ must also be contained in $G$. For instance, let these be denoted by $\mathfrak{R}_{i_{1}}$; let $\mathfrak{R}_{i_{1} i_{2}}$ be not coprime to these; in general, let $\mathfrak{R}_{i_{1} \dots i_{\lambda}}$ be not coprime to $\mathfrak{R}_{i_{1} \dots i_{\lambda - 1}}$. Because we are only dealing with finitely many ideals in total, this procedure must terminate in finitely many steps, that is, there comes a point after which no ideals different from all those preceding are yielded; the thus obtained system of ideals contructs a group $G$ with the desired properties, because all ideals not contained in $G$ are by construction coprime to all those contained in $G$. Suppose there further exists a splitting into two subgroups $G^{(1)}$ and $G^{(2)}$; and without loss of generality let $\mathfrak{R}_{i_{1} \dots i_{\lambda}}$ be an element of $G^{(1)}$. Because the relationship of coprimality is symmetric, $G^{(1)}$ must then also contain $\mathfrak{R}_{i_{1} \dots i_{\lambda - 1}}, \dots , \mathfrak{R}_{i_{1}}$, and so also $\mathfrak{R}$ and consequently the whole group $G$, which proves irreducibility. If we then proceed accordingly with the groups not contained in $G$, then we obtain a grouping $G_{1}, \dots, G_{\tau}$ of all $\mathfrak{R}$. This grouping is \emph{unique}; let $G_{1}^{'}, \dots, G_{\tau}^{'}$ be a second grouping and $\mathfrak{R}_{i_{1} \dots i_{\lambda}}$ an element of $G_{i}^{'}$. Then $G_{i}^{'}$ also contains $\mathfrak{R}$ and hence $G_{1}$, and because of the irreducibility requirement cannot contain any elements different from $G_{1}$; thus $G_{i}^{'}$ is the same as $G_{1}$.

Now let $\mathfrak{T}_{i}$ be the least common multiple of the ideals $\mathfrak{R}$ incorporated into a group $G_{i}$. We show that $\mathfrak{M} = [\mathfrak{T}_{1}, \dots, \mathfrak{T}_{\tau}]$ is \emph{a representation of $\mathfrak{M}$ through coprime irreducible ideals}. First of all, the decomposition of the $\mathfrak{T}$ into the $\mathfrak{R}$ shows that $[\mathfrak{T}_{1}, \dots, \mathfrak{T}_{\tau}]$ really does give a represention of $\mathfrak{M}$. Furthermore, by Theorem XIV the $\mathfrak{T}$ are pairwise coprime, and each $\mathfrak{T}$ is coprime irreducible. Hence the \emph{existence} of such a representation is proven.

For the proof of uniqueness, let $\mathfrak{M} = [\bar{\mathfrak{T}}_{1}, \dots, \bar{\mathfrak{T}}_{\bar{\tau}}]$ be a second such representation. If we decompose the $\bar{\mathfrak{T}}$ into their relatively prime irreducible ideals $\mathfrak{R}$, then the $\mathfrak{R}$ appearing in different $\bar{\mathfrak{T}}$ are coprime to each other by Theorem XIV, and so are also mutually relatively prime; they are therefore the same as the uniquely defined relatively prime irreducible ideals $\mathfrak{R}$ of $\mathfrak{M}$. By Theorem XIV, the $\bar{\mathfrak{T}}_{i}$ further generate a grouping $G_{i}^{'}$ of the $\mathfrak{R}$ with the given properties. Because, however, this grouping is \emph{unique}, and every $\bar{\mathfrak{T}}_{i}$ is \emph{uniquely} defined by the group $G_{i}^{'} = G_{i}$, it follows that $\bar{\mathfrak{T}}_{i} = \mathfrak{T}_{i}$, \emph{proving uniqueness}.

Furthermore, for pairwise coprime ideals \emph{the least common multiple is the same as the product}.

Because by Theorem XIV the complement $\mathfrak{L}_{i}$ is also coprime to $\mathfrak{T}_{i}$ for $\mathfrak{M} = [\mathfrak{T}_{1} \dots \mathfrak{T}_{\tau}]$, there therefore exist elements
\begin{align*}
t_{i} \equiv 0\ (\mathfrak{T}_{i}); \quad l_{i} \equiv 0\ (\mathfrak{L}_{i}); \quad \varepsilon = t_{i} + l_{i}.
\end{align*}
It therefore follows from $f \equiv 0\ (\mathfrak{T}_{i}), f \equiv 0\ (\mathfrak{L}_{i})$ that, because
\begin{align*}
f = f \varepsilon = f t_{i} + f l_{i},\ \text{it also holds that}\ f \equiv 0\ (\mathfrak{L}_{i} \cdot \mathfrak{T}_{i}).
\end{align*}
Because conversely $\mathfrak{L}_{i} \cdot \mathfrak{T}_{i}$ is divisible by $[\mathfrak{L}_{i}, \mathfrak{T}_{i}]$, it follows that $[\mathfrak{L}_{i}, \mathfrak{T}_{i}] = \mathfrak{L}_{i} \cdot \mathfrak{T}_{i}$, and by extension of the procedure over the $\mathfrak{L}_{i}$, we conclude that
\begin{align*}
\mathfrak{M} = [\mathfrak{T}_{1}, \dots, \mathfrak{T}_{\tau}] = \mathfrak{T}_{1} \cdot \mathfrak{T}_{2} \cdot \dots{} \cdot \mathfrak{T}_{\tau}.
\end{align*}
We have therefore proved \\

\textbf{Theorem XV}. \emph{Every ideal can be uniquely expressed as the product of finitely many pairwise coprime and coprime irreducible ideals}.

\section*{\S 9. Development of the study of modules. Equality of the number of components in decompositions into irreducible modules.}

\tab We now show that the content of the first three sections, which relates to \emph{irreducible} ideals, not primary and prime ideals, still holds under less restrictive conditions. These sections in particular do not use the law of commutativity of multiplication and concern only the property of ideals being modules, and so remain upheld for modules over non-commutative rings, which are now to be defined. The definition of these modules shall be based upon a \emph{\underline{double domain}} $(\Sigma, T)$ with the following properties:

$\Sigma$ \emph{is an abstractly defined non-commutative ring}, that is, $\Sigma$ is a system of elements $a, b, c, \dots$, for which two operations are defined; ring addition ($\hash$) and ring multiplication ($\bigtimes$), which satisfy the laws set out in \S 1, with the exception of law 4. regarding the commutativity of ring multiplication.

$T$ is a system of elements $\alpha, \beta, \gamma, \dots$ for which, in conjunction with $\Sigma$, two operations are also defined; \emph{addition}, which for every two elements $\alpha, \beta$ uniquely generates a third $\alpha + \beta$, and \emph{multiplication of an element $\alpha$ of $T$ with an element $c$ of $\Sigma$}, which uniquely generates an element $c \cdot \alpha$ of $T$.\footnote{We are therefore dealing with \textquotedblleft right" multiplication, a \textquotedblleft right" domain $T$, and thus \textquotedblleft right" modules and ideals. Were we to have based $T$ on a left multiplication $\alpha \cdot c$, then a corresponding theory of the left modules and ideals would follow; $M$ would contain in addition to $\alpha$ also $\alpha \cdot c$. The law of associativity would have the form $(\gamma \cdot b) \cdot a = \gamma \cdot (b \bigtimes a)$ here.}

The following laws apply to these operations: \\
1. The law of associativity of addition: $(\alpha + \beta) + \gamma = \alpha + (\beta + \gamma)$. \\
2. The law of commutativity of addition: $\alpha + \beta = \beta + \alpha$. \\
3. The law of unrestricted and unique subtraction: there exists one and only one element $\xi$ in $T$ which satisfies the equation $\alpha + \xi = \beta$ (written $\xi = \beta - \alpha$). \\
4. The law of associativity of multiplication: $a \cdot (b \cdot \gamma) = (a \bigtimes b) \cdot \gamma$. \\
5. The law of distributivity: $(a \hash b) \cdot \gamma = a \cdot \gamma + b \cdot \gamma; \quad c \cdot (\alpha + \beta) = c \cdot \alpha + c \cdot \beta$.

The existence of the zero element follows from these conditions, as is well-known, and also the validity of the law of distributivity for subtraction and multiplication:
\begin{align*}
(a \dotdiv b) \cdot \gamma = a \cdot \gamma - b \cdot \gamma; \quad c \cdot (\alpha - \beta) = c \cdot \alpha - c \cdot \beta;
\end{align*}
where ($\dotdiv$) denotes subtraction in $\Sigma$. If $\Sigma$ contains a unit $\varepsilon$, then $\varepsilon \cdot \alpha = \alpha$ holds for all elements $\alpha$ of $T$.

Let a \emph{module} $M$ over $(\Sigma, T)$ be understood to be a system of elements of $T$ which satisfies the following two conditions: \\
1. \emph{For each element $\alpha$ of $M$, $c \cdot \alpha$ is also an element of $M$, where $c$ is an arbitrary element of $\Sigma$.} \\
2. \emph{For each pair of elements $\alpha$ and $\beta$ of $M$, the difference $\alpha - \beta$ is also an element of $M$}; therefore for each element $\alpha$ of $M$, $n \alpha$ is also an element of $M$ for every integer $n$.\footnote{These integers are again to be considered as abbreviatory symbols, not as ring elements.} 

By this definition, \emph{$T$ itself constitutes a module in $(\Sigma, T)$}. If in particular $T$ and the operations defined for it coincide with the ring $\Sigma$ and the operations applicable to it, then the module $M$ becomes a (right) ideal $\mathfrak{M}$ in $\Sigma$. If $\Sigma$ is taken to be commutative, then the usual concept of an ideal arises, which therefore comes about as a special case of the concept of a module.\footnote{The simplest example of a module is the module consisting of \underline{integer linear forms}; here $\Sigma$ consists of all integers and $T$ consists of all integer linear forms. A somewhat more general module arises if algebraic integers are used instead of integers in $\Sigma$ and $T$, or for instance all even numbers. If we consider the complex of all coefficients as \emph{one} element each time instead of the linear forms, then the operations in $\Sigma$ and $T$ are in fact different. Ideals in non-commutative rings of polynomials form the subject matter of the joint work of Noether-Schmeidler. The lectures of Hurwitz on the number theory of quaternions (Berlin, Springer 1919) and the works of Du Pasquier cited there relate to ideals in further special non-commutative rings.}

All definitions in \S 1 remain upheld for modules: so $\alpha \equiv 0\ (M)$ and $N \equiv 0\ (M)$ mean that $\alpha$ and each element of $N$ respectively are in $M$; in other words, $\alpha$ and $N$ respectively are divisible by $M$. $M$ is a proper divisor of $N$ if $M$ contains elements not in $N$; it follows from $N \equiv 0\ (M)$, $M \equiv 0\ (N)$ that $M = N$. The definitions of the greatest common divisor and the least common multiple remain upheld word for word. In particular, if $M$ contains a finite number of elements $\alpha_{1} \dots \alpha_{\varrho}$ such that $M = (\alpha_{1} \dots \alpha{\varrho})$, that is $\alpha = c_{1} \alpha_{1} + \dots + c_{\varrho} \alpha_{\varrho} + n_{1} \alpha_{1} + \dots n_{\varrho} \alpha_{\varrho}$ for every $\alpha \equiv 0\ (M)$, where the $c_{i}$ are elements of $\Sigma$ and the $n_{i}$ are integers, then $M$ is called a \emph{finite module}, and $\alpha_{1} \dots \alpha_{\varrho}$ a \emph{module basis}.

In the following we now use, analogously to \S 1, only \emph{domains $(\Sigma, T)$ that satisfy the finiteness condition: every module in $(\Sigma, T)$ is finite, and so has a module basis}.

Theorem I of the Finite Chain then also holds for modules for this domain $(\Sigma, T)$, as the proof there shows, and thus all requirements for \S \S 2 and 3 are satisfied. It remains to directly carry over Definition I and Lemma I regarding shortest and reduced representations, and likewise Theorem II regarding the ability to represent each module as the least common multiple of finitely many irreducible modules, whereby Lemma II shows that each such shortest representation is immediately reduced. Furthermore, Theorem III also remains upheld, which conveys the reducibility of a module through the properties of its complement; and from that follows Lemma III, and finally Theorem IV, which states the \emph{equality of the number of components in two different shortest representations of a module as the least common multiple of irreducible modules}. Theorem IV then yields Lemma IV as the converse of Lemma I on reduced representations. \\

\textbf{Remark}. The same reasoning shows that all of these theorems and definitions remain upheld if we understand all modules to be \emph{two-sided} (that is, if $\alpha$ is an element, then $c \cdot \alpha$ and $\alpha \cdot c$ are also elements, and if $\alpha$ and $\beta$ are elements, then $\alpha - \beta$ is also an element) for two-sided \underline{domains} $T$, that is, domains that are both left domains and right domains. \\

While all these theorems are based only on the notions of \emph{divisibility} and the \emph{least common multiple}, the further theorems of uniqueness are based crucially on the \emph{concept of product}, and therefore do not have a direct translation. Therefore the definitions of primary and prime ideals do not translate to modules, because the product of two elements of $T$ is not defined. Although it is feasible to formally construct a translation to non-commutative rings, it loses its meaning, because here the existence of the associated prime ideals cannot be proved,\footnote{For example, from
\begin{align*}
p_{1}^{\lambda_{1}} \equiv 0\ (\mathfrak{Q}), \quad p_{2}^{\lambda_{2}} \equiv 0\ (\mathfrak{Q})\quad \text{it does not follow that}\quad (p_{1} - p_{2})^{\lambda_{1} + \lambda_{2}} \equiv 0\ (\mathfrak{Q})
\end{align*}
and likewise from
\begin{align*}
(a \cdot b)^{\lambda} \equiv 0\ (\mathfrak{Q})\quad \text{it does not follow that}\quad a^{\lambda} \cdot b^{\lambda} \equiv 0\ (\mathfrak{Q}).
\end{align*}
As a result, $\mathfrak{P}$ cannot be proven to be an ideal, nor does it have the properties of prime ideals. For two-sided ideals, $\mathfrak{P}$ can be proven to be an ideal, but still not a prime ideal.} and the proof that an irreducible ideal is primary also breaks down. However, these two cases lay the foundations for the uniqueness theorems \--- in contrast, if the non-commutative ring contains a unit, coprime and coprime irreducible ideals can be defined, and the reasoning used in the proof of Theorem II shows that every ideal can be represented as the least common multiple of \emph{finitely many} pairwise coprime and coprime irreducible ideals.\footnote{For special \textquotedblleft completely reduced" ideals, uniqueness theorems can also be established here; cf. the joint work of Noether-Schmeidler.}

Finally, we mention another criterion sufficient for the finiteness condition to be satisfied in $(\Sigma, T)$: \emph{if $\Sigma$ contains a unit and satisfies the finiteness condition, and if $T$ is itself a finite module in $(\Sigma, T)$, then every module in $(\Sigma, T)$ is finite}.

It follows from the existence of the unit in $\Sigma$ namely that for
\begin{align*}
\mathfrak{M} = (f_{1}, \dots, f_{\varrho}),\quad f = \bar{b}_{1} f_{1} + \dots + \bar{b}_{\varrho} f_{\varrho} + n_{1} f_{1} + \dots + n_{\varrho} f_{\varrho},
\end{align*}
where the $\bar{b}_{i}$ are elements of $\Sigma$ and the $n_{i}$ are integers, there is also a representation of the form $f = b_{1} f_{1} + \dots + b_{\varrho} f_{\varrho}$, where the $b_{i}$ are elements of $\Sigma$. This is because $f_{i} = \varepsilon f_{i}$ implies $n_{i} f_{i} = n_{i} \varepsilon f_{i}$, and since $n_{i} \varepsilon = (\varepsilon + \dots + \varepsilon)$ belongs to $\Sigma$, it then follows that $b_{i} = \bar{b}_{i} + n_{i} \varepsilon$ is an element of $\Sigma$. The requirement regarding $T$ now states that every element $\alpha$ of $T$ has a representation
\begin{align*}
\alpha &= \bar{a}_{1} \tau_{1} + \dots + \bar{a}_{k} \tau_{k} + n_{1} \tau_{1} + \dots + n_{k} \tau_{k}, \\
\text{and hence}\quad \quad \alpha &= a_{1} \tau_{1} + \dots + a_{k} \tau_{k},
\end{align*}
where the second representation arises from the first because $\varepsilon \tau_{i} = \tau_{i}$, as above.

If the $\alpha$ now run through the elements of a module $M$ in $(\Sigma, T)$, then the coefficients $a_{k}$ of $\tau_{k}$ run through the elements of an ideal $\mathfrak{M}_{k}$ in $\Sigma$; by the above, for each $a_{k} \equiv 0\ (\mathfrak{M}_{k})$, we also have $a_{k} = b_{1} a_{k}^{(1)} + \dots + b_{\varrho} a_{k}^{(\varrho)}$. If we let $\alpha^{(i)}$ denote an element of $M$ for which the coefficient of $\tau_{k}$ is equal to $a_{k}^{(i)}$, then $\alpha - b_{1} \alpha^{(1)} - \dots{} - b_{\varrho} \alpha^{(\varrho)}$ is an element belonging to $M$ which depends only on $\tau_{1}, \dots, \tau_{k-1}$. This procedure can be repeated on the collection of these elements no longer containing $\tau_{k}$, which construct a module $M^{'}$, so that after finitely many repetitions \emph{the proposition is proved}.

\section*{\S 10. Special case of the polynomial ring.}

\tab \textbf{1}. The ring $\Sigma$ that we take as our starting point consists of all polynomials in $x_{1}, \dots , x_{n}$ with arbitrary complex coefficients, for which the finiteness condition is satisfied due to Hilbert's Module Basis Theorem (Math. Ann., vol. 36). This section concerns the \emph{connection of our theorems with the known theorems of elimination theory and module theory}.

This connection is established through the following special case of a famous theorem of Hilbert's:\footnote{\"{U}ber die vollen Invariantensysteme. Math. Ann. 42 (1893), \S 3, p313.}

If $f$ vanishes for every (finite) system of values of $x_{1}, \dots, x_{n}$ that is a root of all polynomials of a prime ideal $\mathfrak{P}$ (we call such a system a root of $\mathfrak{P}$), then $f$ is divisible by $\mathfrak{P}$. In other words, a prime ideal $\mathfrak{P}$ consists of \emph{all} polynomials that vanish at these roots.\footnote{This special case can also be proven directly in the case of homogeneous forms, as Lasker [Math. Ann. 60 (1905), p607] has shown, and then conversely this again implies the Hilbert Theorem (in the homogeneous and the inhomogeneous case), which we can state as follows: if an ideal $\mathfrak{R}$ vanishes at all roots of $\mathfrak{M}$, then a power of $\mathfrak{R}$ is divisible by $\mathfrak{M}$. This theorem, and likewise the special case, only holds however if the codomain of the $x$ is \emph{algebraically closed}, and therefore cannot follow from our theorems alone, but must use the existence of roots; for instance, in the case that an ideal for which the only root is $x_{1} = 0, \dots, x_{n} = 0$ contains all products of powers of the $x$ of a particular dimension. The remaining proof can be simplified somewhat by using our theorems via Lasker. Lasker must in particular also make use of the Hilbert Theorem for the proof of the decomposition of an ideal into maximal primary ideals.} 

If a product $f \cdot g$ vanishes for all roots of $\mathfrak{P}$, then at least one factor vanishes; the roots form an \emph{\underline{irreducible algebraic figure}}. Conversely, should we begin with this definition of the irreducible figure, then it follows that the collection of polynomials vanishing on an irreducible figure forms a prime ideal; prime ideals and irreducible figures therefore correspond with each other bijectively. Furthermore, because $\mathfrak{Q} \equiv 0\ (\mathfrak{P})$ and $\mathfrak{P}^{\varrho} \equiv 0\ (\mathfrak{Q})$, the roots of a primary ideal are the same as those of its associated prime ideal.\footnote{Macaulay (cf. Introduction) uses this property of a primary ideal having an irreducible figure in the definition, while Lasker only incorporates the concept of the manifold of a figure in the definition, which is otherwise abstractly defined. The primary ideals which vanish only for $x_{1} = 0, \dots, x_{n} = 0$ take a special place in Lasker.} The representation of an ideal as the least common multiple of maximal primary ideals therefore yields a partition of all roots of the ideal into irreducible figures; and as Lasker has shown, the converse also holds. \emph{The proof of the uniqueness of the associated prime ideals therefore corresponds here with the Fundamental Theorem of Elimination Theory regarding the unique decomposability of an algebraic figure into irreducible figures}; it can serve as the equivalent of this theorem of elimination theory for special polynomial rings where no unique representation of the polynomial as the product of irreducible polynomials of the ring exists, and consequently also no elimination theory.

The irreducible figures corresponding to the isolated primary ideals are exactly those occurring in the \textquotedblleft minimal resolvent";\footnote{Cf. for example J. K\"{o}nig, Einleitung in die allgemeine Theorie der algebraischen Gr\"{o}\ss en (Leipzig, Teubner, 1903), p235.} since the roots of each non-isolated maximal primary ideal are at the same time roots of at least one isolated one, namely one for which the associated prime ideal is divisible by that of the non-isolated ideal. The uniqueness of the isolated primary ideals therefore yields new invariant multiplicities in the exponents. Additionally, the uniqueness theorems for the decomposition of the primary ideals into irreducible ideals can be viewed as an addendum to elimination theory, in the sense of multiplicity.

Following these remarks, the different decompositions can be interpreted in their relation to the algebraic figures. The pairwise coprime ideals correspond to figures that have no roots in common; for the mutually relatively prime ideals, likewise no irreducible figure of one ideal has roots in common with that of another; the maximal primary ideals vanish only in irreducible figures which are all different from each other; for the decomposition into irreducible ideals, the same irreducible figures can also appear repeatedly. 

It should be noted that the ring of all homogeneous forms can also be used in place of the general polynomial ring, since it is easy to convince ourselves that the general theorems also remain upheld for the operations valid there\footnote{The following example shows that here inhomogeneous decompositions can also exist as well as homogeneous ones in the case of non-uniqueness, however:
\begin{align*}
(x^{3}, xy, y^{3}) = [(x^{3}, y); (y^{3}, x)] = [(xy, x^{3}, y^{3}, x + y^{2}); (xy, x^{3}, y^{3}, y + x^{2})].
\end{align*}} \--- the addition is only defined for forms of the same dimension.

A simple example of the four different decompositions \--- for which the formulae below follow \--- is given, following the above remarks, by one straight line and two further straight lines skew to it and intersecting each other, one of which contains a point of higher multiplicity other than the point of intersection. The decomposition into coprime-irreducible ideals corresponds to the decomposition into the straight line and the figure skew to it; this figure splits into the two straight lines it is composed of for the decomposition into relatively prime irreducible ideals; the decomposition into maximal primary ideals corresponds to a detachment of the point of higher multiplicity, while the decomposition into irreducible ideals requires the removal of this point. 

Should we take this point as the starting point, and the straight line passing through it as the $y$-axis, the straight line intersecting this parallel to the $x$-axis, and the straight line skew to it parallel to the $z$-axis, then one such configuration is represented by the following \emph{irreducible ideals}:\footnote{The first three of these are irreducible since they are prime ideals; $\mathfrak{B}_{4}$ is irreducible because it only has the divisors $(x^{2}, y, z)$ and $(x, y, z)$; $\mathfrak{B}_{5}$ is irreducible because \emph{every} divisor contains the polynomial $xy$.}
\begin{equation*}
\mathfrak{B}_{1} = (x - 1, y); \quad \mathfrak{B}_{2} = (y - 1, z); \quad \mathfrak{B}_{3} = (x, z); \quad \mathfrak{B}_{4} = (x^{3}, y, z); 
\end{equation*}
\begin{equation*}
\mathfrak{B}_{5} = (x^{2}, y^{2}, z).
\end{equation*}
The \emph{associated prime ideals} are:
\begin{align*}
\mathfrak{P}_{1} = \mathfrak{B}_{1}; \quad \mathfrak{P}_{2} = \mathfrak{B}_{2}; \quad \mathfrak{P}_{3} = \mathfrak{B}_{3}; \quad \mathfrak{P}_{4} = \mathfrak{P}_{5} = (x, y, z).
\end{align*}
The \emph{maximal primary ideals} are:
\begin{align*}
\mathfrak{Q}_{1} = \mathfrak{B}_{1}; \quad \mathfrak{Q}_{2} = \mathfrak{B}_{2}; \quad \mathfrak{Q}_{3} = \mathfrak{B}_{3}; \quad \mathfrak{Q}_{4} = [\mathfrak{B}_{4}, \mathfrak{B}_{5}] = (x^{3}, y^{2}, x^{2}y, z).
\end{align*}
The \emph{relatively prime irreducible ideals} are:
\begin{align*}
\mathfrak{R}_{1} = \mathfrak{Q}_{1}; \quad \mathfrak{R}_{2} = \mathfrak{Q}_{2}; \quad \mathfrak{R}_{3} = [\mathfrak{Q}_{3}, \mathfrak{Q}_{4}] = (x^{3}, x^{2}y, xy^{2}, z).
\end{align*}
The \emph{coprime irreducible ideals} are:
\begin{align*}
\mathfrak{S}_{1} = \mathfrak{R}_{1}; \quad \mathfrak{S}_{2} = [\mathfrak{R}_{2}, \mathfrak{R}_{3}] = ((y - 1)x^{3}, (y - 1)x^{2}y, (y - 1)xy^{2}, z).
\end{align*}
This produces the \emph{total ideal}: 
\begin{align*}
\mathfrak{M} &= [\mathfrak{S}_{1}, \mathfrak{S}_{2}] \\ 
&= \mathfrak{S}_{1} \cdot \mathfrak{S}_{2} \\
&= ((x - 1)(y - 1)x^{3}, (y - 1)x^{2}y, (y - 1)xy^{2}, (x - 1)z, y(y-1)x^{3}, yz),
\end{align*}
which gives
\begin{equation*}
1 = -(y - 1)x^{3} + (y-1)(x^{3} - 1) + y, \quad (y - 1)(x^{3} - 1) + y \equiv 0\ (\mathfrak{S}_{1});
\end{equation*}
\begin{equation*}
 -(y - 1)x^{3} \equiv 0\ (\mathfrak{S}_{2}).
\end{equation*}

Here the ideals $\mathfrak{B}_{1}, \mathfrak{B}_{2}, \mathfrak{B}_{3}$ are isolated, and so also uniquely determined in the decompositions into irreducible and maximal primary ideals. The ideals $\mathfrak{B}_{4}$ and $\mathfrak{B}_{5}$, and respectively $\mathfrak{Q}_{4}$, are non-isolated; they are not uniquely determined, but rather can be replaced for example with $\mathfrak{D}_{4} = (x^{3}, y + \lambda x^{2}, z),\ \mathfrak{D}_{5} = (x^{2} + \mu xy, y^{2}, z)$; similarly, $\mathfrak{Q}_{4}$ can be replaced with
\begin{align*}
\bar{\mathfrak{Q}}_{4} = (x^{3}, x^{2}y, y^{2} + \lambda xy, z).
\end{align*}

\textbf{2}. Similarly to the general (and the integer) polynomial ring, every \emph{finite integral domain of polynomials} also satisfies the finiteness condition\footnote{Conversely, if the finiteness condition is satisfied for a polynomial ring, and if each polynomial has at least one representation where the factors of lower degree are in $x$, then the ring is a finite integral domain.} \--- as Hilbert's Module Basis Theorem shows \--- where the coefficients can be assigned an arbitrary field. We shall give another example for the \emph{ring of all even polynomials}, as the simplest ring where, because $x^{2} \cdot y^{2} = (xy)^{2}$, no unique factorisation of the polynomial into irreducible polynomials of the ring exists. It uses the same configuration as in the above example, which is now given through the \emph{irreducible ideals}
\begin{align*}
\mathfrak{B}_{1} &= (x^{2} - 1, xy, y^{2}, yz); \quad \mathfrak{B}_{2} = (y^{2} - 1, xz, yz, z^{2}); \\
\mathfrak{B}_{3} &= (x^{2}, xy, xz, yz, z^{2}); \quad \mathfrak{B}_{4} = (x^{4}, xy, y^{2}, xz, yz, z^{2}); \\
\mathfrak{B}_{5} &= (x^{2}, y^{2}, xz, yz, z^{2}).
\end{align*}

The \emph{associated prime ideals} are:
\begin{align*}
\mathfrak{P}_{1} = \mathfrak{B}_{1}; \quad \mathfrak{P}_{2} = \mathfrak{B}_{2}; \quad \mathfrak{P}_{3} = \mathfrak{B}_{3}; \quad \mathfrak{P}_{4} = \mathfrak{P}_{5} = (x^{2}, xy, y^{2}, xz, yz, z^{2}).
\end{align*}

It follows that $\mathfrak{P}_{1}$ is a prime ideal because every polynomial of the ring has the following form:
\begin{align*}
f \equiv \phi (z^{2}) + xz \psi (z^{2})\ (\mathfrak{P}_{1}).
\end{align*}
Therefore let
\begin{align*}
f_{1} \equiv \phi_{1} (z^{2}) + xz \psi_{1} (z^{2}); \quad f_{2} \equiv \phi_{2} (z^{2}) + xz \psi_{2} (z^{2}),
\end{align*}
and thus, since $f_{1} \cdot f_{2} \equiv 0\ (\mathfrak{P}_{1})$, we have the existence of the following equations:
\begin{align*}
\phi_{1} \phi_{2} + z^{2} \psi_{1} \psi_{2} = 0; \quad \phi_{2} \psi_{1} + \phi_{1} \psi_{2} = 0,
\end{align*}
and so $f_{1} \equiv 0\ (\mathfrak{P}_{1})$ or $f_{2} \equiv 0\ (\mathfrak{P}_{1})$.

We can show that $\mathfrak{P}_{2}$ is a prime ideal in precisely the same way; $\mathfrak{P}_{3}$ is a prime ideal because every polynomial of the ring mod $\mathfrak{P}_{3}$ is congruent to a polynomial in $y^{2}$; $\mathfrak{P}_{4}$ consists of all polynomials of the ring. It therefore also follows that $\mathfrak{B}_{1}$, $\mathfrak{B}_{2}$ and $\mathfrak{B}_{3}$ are irreducible, as they are prime ideals; $\mathfrak{B}_{4}$ and $\mathfrak{B}_{5}$ each have only the sole proper divisor $\mathfrak{P}_{4}$, and so are necessarily irreducible.

From the irreducible ideals arise the \emph{maximal primary ideals}:
\begin{align*}
\mathfrak{Q}_{1} = \mathfrak{B}_{1}; \ \ \mathfrak{Q}_{2} = \mathfrak{B}_{2}; \ \ \mathfrak{Q}_{3} = \mathfrak{B}_{3}; \ \ \mathfrak{Q}_{4} = [\mathfrak{B}_{4}, \mathfrak{B}_{5}] = (x^{4}, x^{3}y, y^{2}, xz, yz, z^{2});
\end{align*}
the \emph{relatively prime irreducible ideals}:
\begin{align*}
\mathfrak{R}_{1} = \mathfrak{Q}_{1}; \quad \mathfrak{R}_{2} = \mathfrak{Q}_{2}; \quad \mathfrak{R}_{3} = [\mathfrak{Q}_{3}, \mathfrak{Q}_{4}] = (x^{4}, x^{3}y, x^{2}y^{2}, xy^{3}, xz, yz, z^{2}),
\end{align*}
the \emph{coprime irreducible ideals}:
\begin{align*}
\mathfrak{S}_{1} &= \mathfrak{R}_{1} \\
\mathfrak{S}_{2} &= [\mathfrak{R}_{2}, \mathfrak{R}_{3}] \\
&= ((y^{2} - 1)x^{4}, (y^{2} - 1)x^{3}y, (y^{2} - 1)x^{2}y^{2}, (y^{2} - 1)xy^{3}, xz, yz, z^{2}).
\end{align*}
As in the first example, it follows that
\begin{align*}
[\mathfrak{B}_{4}, \mathfrak{B}_{5}] = [\mathfrak{D}_{4}, \mathfrak{D}_{5}];
\end{align*}
where $\mathfrak{D}_{4} = (x^{4}, xy + \lambda x^{2}, \dots )$ and $\mathfrak{D}_{5} = (x^{2} + \mu xy, \dots )$ for $\lambda \cdot \mu \neq 1$; and
\begin{align*}
[\mathfrak{Q}_{3}, \mathfrak{Q}_{4}] = [\mathfrak{Q}_{3}, \bar{\mathfrak{Q}}_{4}]; \quad \text{where}\ \bar{\mathfrak{Q}}_{4} = (x^{4}, x^{3}y, y^{2} + \lambda xy, \dots ).
\end{align*}

The remaining irreducible ideals (and maximal primary ideals respectively) $\mathfrak{B}_{1}$, $\mathfrak{B}_{2}$, $\mathfrak{B}_{3}$ are uniquely determined as isolated ideals.

\section*{\S 11. Examples from number theory and the theory of differential expressions.}

\tab \textbf{1}. Let the ring $\Sigma$ consist of all \emph{even integers}. $\Sigma$ can then be bijectively assigned to all of the integers, since each number $2a$ in $\Sigma$ can be allocated the number $a$. From this it immediately follows that every ideal in $\Sigma$ is a principal ideal $(2a)$, where in the basis representation $2c = n \cdot 2a$ of each element $2c$ of the ideal, the odd numbers $n$ only amount to abbreviations for finite sums.

The \emph{prime ideals} of the ring are given by $\mathfrak{P}_{0} = \mathfrak{O} = (2)$ and $\mathfrak{P} = (2p)$, where $p$ is an odd prime number; therefore every prime ideal is divisible by $\mathfrak{P}_{0}$, but by no other prime ideal. The \emph{primary ideals} are given by $\mathfrak{Q}_{\varrho_{0}} = (2 \cdot 2^{\varrho_{0}})$ and $\mathfrak{Q}_{\varrho} = (2p^{\varrho})$; they are at the same time \emph{irreducible ideals}, and by what has been said about prime ideals, each two corresponding to different odd prime numbers are mutually relatively prime, but no $\mathfrak{Q}$ is relatively prime to any $\mathfrak{Q}_{\varrho_{0}}$,\footnote{In fact, it always follows from $2b \cdot 2p_{1}^{\varrho_{1}} \equiv 0\ (2p_{2}^{\varrho_{2}})$ for odd $p_{1} \neq p_{2}$ that $2b \equiv 0\ (2p_{2}^{\varrho_{2}})$; however, $2b \cdot 2p^{\varrho} \equiv 0\ (2 \cdot 2^{\varrho_{0}})$ only implies that $2b \equiv 2 \cdot 2^{\varrho_{0} - 1}\ (2 \cdot 2^{\varrho_{0}})$.} and so the $\mathfrak{Q}_{\varrho_{0}}$ are the only non-isolated primary ideals. The unique decomposition of $a$ into prime powers corresponds to the \emph{unique} representation of the ideal $(2a)$ through maximal primary (and at the same time irreducible) ideals:
\begin{align*}
(2a) = [(2 \cdot 2^{\varrho_{0}}), (2p_{1})^{\varrho_{1}}, \dots, (2p_{\alpha})^{\varrho_{\alpha}}];
\end{align*}
therefore, contrary to the examples from the polynomial ring, the \emph{non-isolated} maximal primary ideals are also \emph{uniquely} determined. As for $\varrho_{0} = 0$, it behaves also as a representation through mutually relatively prime ideals, while for $\varrho_{0} > 0$ the ideal is relatively prime irreducible.

While the four different decompositions therefore coincide in the ring of all integers, this is only the case here for the two decompositions into maximal primary ideals and irreducible ideals on the one hand, and for the decompositions into coprime irreducible ideals (every ideal is coprime irreducible, because the ring has no units) and relatively prime irreducible ideals for $\varrho_{0} > 0$ on the other hand, while for $\varrho_{0} = 0$ the coprime irreducible and relatively prime irreducible decompositions are different from each other. At the same time, an example presents itself here. A prime ideal can be divisible by another without having to be identical to it; more generally, \emph{the factorisation does not follow from divisibility}. The last one \--- as a consequence of the fact that the ring contains no units \--- is also the reason that no unique factorisation of the numbers of the ring into irreducible numbers of the ring exists, although each ideal is a principal ideal; the introduction of the least common multiple therefore proves to be necessary here. It should also be noted that the relationships remain exactly \emph{the same} if instead of all even numbers, we use \emph{all numbers divisible by a fixed prime number or prime power}.

However, \emph{irreducible and primary ideals are different} as well if $\Sigma$ consists of \emph{all numbers divisible by a composite number $g = \mathfrak{p}_{1}^{\sigma_{1}} \dots \mathfrak{p}_{\nu}^{\sigma_{\nu}}$}. Every ideal is again a principal ideal $(g \cdot a)$, and the \emph{prime ideals} are again given by $\mathfrak{P}_{0} = \mathfrak{O} = (g)$, $\mathfrak{P} = (g \cdot p)$, where $p$ is a prime number different from the prime numbers $\mathfrak{p}$ appearing in $g$. In contrast, the \emph{irreducible ideals} are given by $\mathfrak{B}_{\lambda_{i}} = (g \cdot \mathfrak{p}_{i}^{\lambda_{i}})$, $\mathfrak{B}_{\varrho} = (g \cdot p^{\varrho})$; the \emph{primary ideals} are given by $\mathfrak{Q}_{\varrho} = \mathfrak{B}_{\varrho}$ and by the $\mathfrak{Q}_{\lambda_{1} \dots \lambda_{\nu}} = (g \cdot \mathfrak{p}_{1}^{\lambda_{1}} \dots \mathfrak{p}_{\nu}^{\lambda_{\nu}})$ different from the irreducible ideals, where the $\mathfrak{B}_{\lambda_{i}}$ and $\mathfrak{Q}_{\lambda_{1} \dots \lambda_{\nu}}$ all have the same associated prime ideal $\mathfrak{P}_{0} = (g)$. Uniqueness of the decomposition into irreducible ideals also holds here, and consequently uniqueness of the decomposition into maximal primary ideals too, and so the non-isolated ideals are again uniquely determined here too.

\textbf{2}. One example of a \emph{non-commutative ring} is presented by the ideal theory in non-commutative polynomial rings discussed in the Noether-Schmeidler paper. It concerns in particular \textquotedblleft completely reducible" ideals, that is, ideals for which the components of the decomposition are pairwise coprime and have no proper divisors; the components are therefore a fortiori irreducible. The \emph{equality of the number of components} in two different decompositions therefore follows from \S 9, in addition to the isomorphism proved there. Thus a consequence of the decomposition of systems of partial or ordinary linear differential expressions arising as a special case of the paper is obtained, which even appears not to have been remarked upon in the well-known case of an ordinary linear differential expression.

Meanwhile, the system $T$ of all cosets of a fixed ideal $\mathfrak{M}$ together with the non-commutative polynomial ring $\Sigma$ gives a double domain $(\Sigma, T)$, where $T$ has the module property with respect to $\Sigma$, since the difference of two cosets is again a coset, and likewise the product of a coset with an arbitrary polynomial, whereas the product of two cosets does not exist (loc. cit. \S 3). The systems of cosets denoted there as \textquotedblleft subgroups" form examples of modules in double domains $(\Sigma, T)$ where the ring $\Sigma$ is non-commutative.

\section*{\S 12. Example from elementary divisor theory.}

\tab This section deals with a concept of elementary divisor theory contingent on the general developments thereof, which is however itself presumed to be \emph{known}.

Let $\Sigma$ be the ring of all integer matrices with $n^{2}$ elements, for which addition and multiplication are defined in the conventional sense for matrices. $\Sigma$ is then a \emph{non-commutative ring}; the ideals are thus in general one-sided, and two-sided ideals only arise as a special case.\footnote{The ideal theory of these rings forms the subject of the papers of Du Pasquier: Zahlentheorie der Tettarionen, Dissertation Z\"{u}rich, Vierteljahrsschr. d. Naturf. Ges. Z\"{u}rich, 51 (1906). Zur Theorie der Tettarionenideale, ibid., 52 (1907). The content of the two papers is the proof that each ideal is a principal ideal.}

We first show that \emph{every ideal} is a \emph{principal ideal}. In order to do this, we allocate (for right ideals) each matrix $A = (a_{ik})$ a module
\begin{align*}
\mathcal{A} = (a_{11}\xi_{1} + \dots + a_{1n}\xi_{n}, \dots , a_{n1}\xi_{1} + \dots + a_{nn}\xi_{n})
\end{align*}
consisting of integer linear forms. Conversely, each matrix which gives a basis of $\mathcal{A}$ corresponds to this module, so in addition to $A$, $UA$ also corresponds to $\mathcal{A}$, where $U$ is unimodular. More generally, a module $\mathcal{B}$ which is a multiple of $\mathcal{A}$ corresponds to the product $PA$. A single linear form from $\mathcal{A}$ is given by a $P$ which contains only one non-zero row.

Now let $A_{1}, A_{2}, \dots , A_{\nu}, \dots$ be all the elements of an ideal $\mathfrak{M}$, with $\mathcal{A}_{1},$ $\mathcal{A}_{2},$ $\dots ,$ $\mathcal{A}_{\nu},$ $\dots$ the modules assigned to them, $\mathcal{A}$ the greatest common divisor of these, and $UA$ the most general matrix assigned to this module $\mathcal{A}$. To every individual linear form from $\mathcal{A}$ then corresponds a matrix $P_{1}A_{i_{1}} + \dots + P_{\sigma}A_{i_{\sigma}}$ by definition of the greatest common divisor, where as above the $P$ contain only one non-zero row. From this it follows that the matrix $A$ corresponding to a basis of $\mathcal{A}$ can also have such a representation, now with general $P$, and thus is an element of $\mathfrak{M}$. Because furthermore every module $\mathcal{A}_{i}$ is divisible by $\mathcal{A}$, every matrix $A_{i}$ is therefore divisible by $A$; \emph{$A$, and in general $UA$, forms a basis of $\mathfrak{M}$}. If we were to deal with left ideals, then we would have to consider the columns of each matrix as the basis of a module accordingly; each ideal is a principal ideal for which $A$ being a basis implies that $AV$ is also a basis, where $V$ is understood to be an arbitrary unimodular matrix.

In the following we now use \emph{two-sided ideals} in the context of elementary divisor theory, which in particular implies that if $A$ belongs to the ideal, then so does $PAQ$. Thus, by that shown above, the most general basis of one such ideal is given by $UAV$, where $U$ and $V$ are unimodular. The basis elements therefore run through a \emph{class of equivalent matrices},\footnote{The right and left classes respectively corresponding to the basis elements of the one-sided ideals.} and a one-to-one relationship exists between the ideal and the class. Consequently, such a relationship also exists between the ideal and the elementary divisor system $(a_{1} | a_{2} | \dots | a_{n})$ of the class, where the $a_{i}$ are known to be non-negative integers, each of which divides the following one. The matrix of the class occurring in the normal form induced by the elementary divisors can therefore be understood as a special basis of the ideal; the divisibility of the elementary divisors follows from that of the ideals, respectively the classes, and vice versa.

Now, however, Du Pasquier has shown loc. cit. \S 11 that for each two-sided ideal the rank is $n$ and all elementary divisors are the same. In order to be able to consider the case of general elementary divisors, we must therefore start not from the ideals, but directly from the two-sided classes (which shall be denoted by capital letters such as $\mathfrak{A}, \mathfrak{B}, \mathfrak{C}, \dots$). A class $\mathfrak{A} = UAV$ is therefore divisible by another class $\mathfrak{B}$ here if $A = PBQ$.

In general, the following statement holds: \emph{the least common multiple (greatest common divisor) of two classes is obtained through the construction of the least common multiple (greatest common divisor) of the corresponding systems of elementary divisors}.\footnote{In the paper Zur Theorie der Moduln, Math. Ann. 52 (1899), p1, E. Steinitz defines the least common multiple (greatest common divisor) of classes through least common multiples and greatest common divisors of the systems of elementary divisors. Independently of this, the least common multiple of classes is found as the \textquotedblleft congruence composition" in H. Brandt, Komposition der bin\"{a}ren quadratischen Formen relativ einer Grundform, J. f. M. 150 (1919), p1.}

Let $(a_{1} | a_{2} | \dots | a_{n})$, $(b_{1} | b_{2} | \dots | b_{n})$, $(c_{1} | c_{2} | \dots | c_{n})$ be the systems of elementary divisors of $\mathfrak{A}$, $\mathfrak{B}$, $\mathfrak{C}^{*}$ respectively, where $c_{i} = [a_{i}, b_{i}]$, and let $\mathfrak{C} = [\mathfrak{A}, \mathfrak{B}]$. Then $\mathfrak{C}^{*}$ is divisible by $\mathfrak{A}$ and $\mathfrak{B}$, and therefore by $\mathfrak{C}$, and conversely the elementary divisors of $\mathfrak{C}$ are divisible by those of $\mathfrak{C}^{*}$, which implies that $\mathfrak{C}$ is divisible by $\mathfrak{C}^{*}$, and thus $\mathfrak{C} = \mathfrak{C}^{*}$. The proof for the greatest common divisors proceeds likewise.

The unique representation of the elementary divisors $a_{i}$ as the least common multiple of prime powers therefore corresponds to a representation of $\mathfrak{A}$ as the least common multiple of classes $\mathfrak{Q}$, the elementary divisors of which are given by the powers of \emph{one} prime number, or in symbols:
\begin{align*}
\mathfrak{Q} \sim (p^{r_{1}} | p^{r_{2}} | \dots | p^{r_{\varrho}} | 0 | \dots | 0); \quad r_{1} \leq r_{2} \leq \dots \leq r_{\varrho}.
\end{align*}
If in particular the rank $\varrho$ is equal to $n$, then we have a decomposition into coprime and coprime irreducible classes here, which, although the ring is non-commutative, is unique.

The classes $\mathfrak{Q}$ can further be decomposed into classes which correspond to the system of elementary divisors:
\begin{equation*}
\mathfrak{B}_{1} \sim (p^{r_{1}} | \dots | p^{r_{1}});\ \ \mathfrak{B}_{2} \sim (1 | p^{r_{2}} | \dots | p^{r_{2}});\ \dots ;\ \ \mathfrak{B}_{\varrho} \sim (1 | \dots | 1 | p^{r_{\varrho}} | \dots | p^{r_{\varrho}});
\end{equation*}
\begin{equation*}
\mathfrak{B}_{\varrho + 1} \sim (1 | \dots | 1 | 0 | \dots | 0),
\end{equation*}
where the number $1$ appears $(\nu - 1)$ times in each $\mathfrak{B}_{\nu}$. If $r_{1} = r_{2} = \dots{}  = r_{\mu} = 0$ here, then $\mathfrak{B}_{1}, \dots, \mathfrak{B}_{\mu}$ are equal to the trivial class and are therefore left out of the decomposition; the same holds for $\mathfrak{B}_{\varrho + 1}$ if $\varrho = n$. If furthermore $r_{\nu} = r_{\nu + 1} = \dots{} = r_{\nu + \lambda}$, then $\mathfrak{B}_{\nu + 1}, \dots, \mathfrak{B}_{\nu + \lambda}$ are proper divisors of $\mathfrak{B}_{\nu}$, and are therefore likewise to be excluded. Denote those remaining by $\mathfrak{B}_{i_{1}}, \dots, \mathfrak{B}_{i_{k}}$, which now provide a shortest representation, so that $\mathfrak{Q} = [\mathfrak{B}_{i_{1}}, \dots, \mathfrak{B}_{i_{k}}]$ \emph{is the unique decomposition of $\mathfrak{Q}$ into irreducible classes}. Namely, let $\mathfrak{B}_{\nu} \sim (1 | \dots | 1 | p^{r_{\nu}} | \dots | p^{r_{\nu}})$ be representable as the least common multiple of $\mathfrak{C} \sim (1 | \dots | 1 | p^{s_{1}} | \dots | p^{s_{\lambda}})$ and $\mathfrak{D} \sim (1 | \dots | 1 | p^{t_{1}} | \dots | p^{t_{\mu}})$: then the number $1$ must appear in the first $(\nu - 1)$ positions of the system of elementary divisors of $\mathfrak{C}$ and $\mathfrak{D}$; in the $\nu$th position, the exponent $s_{1}$ or $t_{1}$ must be equal to $r_{\nu}$; let this be $s_{1}$ without loss of generality. Because however $r_{\nu} = s_{1} \leq s_{2} \leq s_{\lambda} \leq r_{\nu}$, it follows that $\mathfrak{C} = \mathfrak{B}_{\nu}$, and so $\mathfrak{B}_{\nu}$ is \emph{irreducible}.\footnote{In contrast, the $\mathfrak{B}$ are reducible in one-sided classes; here no unique relationship exists between elementary divisors and classes anymore, and therefore there is also no unique decomposition into irreducible one-sided classes. This is shown by the following example given to me by H. Brandt regarding the decomposition into right classes (where the classes are represented by a basis matrix, and respectively by the corresponding module):
\begin{equation*}
\mathfrak{B} = [\mathfrak{C}_{1}, \mathfrak{C}_{2}] = [\mathfrak{D}_{1}, \mathfrak{D}_{2}]; \quad
\mathfrak{B} \sim \left( \begin{array}{cc}
p & 0 \\
0 & p \\ \end{array} \right) ; \quad
\mathfrak{C}_{1} \sim \left( \begin{array}{cc}
1 & 0 \\
0 & p \\ \end{array} \right) ;
\end{equation*}
\begin{equation*}
\mathfrak{C}_{2} \sim \left( \begin{array}{cc}
p & 0 \\
0 & 1 \\ \end{array} \right) ; \quad
\mathfrak{D}_{1} \sim \left( \begin{array}{cc}
1 & 0 \\
0 & p \\ \end{array} \right) ; \quad
\mathfrak{D}_{2} \sim \left( \begin{array}{cc}
p & 0 \\
(p-1) & 1 \\ \end{array} \right) .
\end{equation*}
In fact, the modules $(\xi, p \eta)$, $(p \xi, \eta)$, $(p \xi, (p-1) \xi + \eta)$ each have only the module $(\xi, \eta)$ as a proper divisor, and so are irreducible and furthermore are different from each other.} The same applies for $\mathfrak{B}_{\varrho + 1}$, where $p$ is replaced with $0$. Each of these irreducible classes gives, as the construction of the least common multiple shows, a particular exponent, and the position where this exponent first appears in the system of elementary divisors of $\mathfrak{Q}$, respectively $\mathfrak{A}$, while $\mathfrak{B}_{\varrho + 1}$ indicates the rank. Because these numbers are uniquely determined by the system of elementary divisors of $\mathfrak{Q}$, and respectively $\mathfrak{A}$, and because the relation between elementary divisors and classes is bijective, the decomposition of $\mathfrak{Q}$, and likewise that of an arbitrary class, into irreducible classes is therefore \emph{unique}. In summary, we can say the following: \emph{Every two-sided class $\mathfrak{A}$ consisting of integer matrices with a bounded number of elements can be uniquely expressed as the least common multiple of finitely many irreducible two-sided classes. Each irreducible class in this represents a fixed prime divisor of the system of elementary divisors of $\mathfrak{A}$, an associated exponent, and the position where this exponent first appears. The irreducible class corresponding to the divisor $0$ indicates the rank of $\mathfrak{A}$}. \\

Erlangen, October 1920. \\
\begin{center}
\footnotesize{(Received on 16/10/1920.)}
\end{center}
\newpage

\section*{Translator's notes.}
\underline{Pairwise coprime} (p2) - This is defined in the usual way. \\
\underline{Relatively prime} (p2) - Although today we consider this to mean the same as pairwise coprimality, Noether distinguishes the two.\\
\underline{Residual module} (p25) - At the time this paper was written, the idea of a quotient was relatively new, and so different names were in use. Lasker's residual module is simply another name for the same concept.\\
\underline{Double domain} (p36) - In modern terms this is a non-commutative ring and a module on it. Noether defines a module on a double domain to be what in modern terms is a submodule of the module in the double domain.\\
\underline{Integer linear forms} (p37) - Also known as integral linear combinations.\\
\underline{Domains} (p38) - This concept does not appear to have a direct modern equivalent. It seems to be similar to the idea of a ring, but more ambiguous in definition. As $T$ is clearly not itself a ring, this name has been kept.\\
\underline{Irreducible algebraic figure} (p40) - This concept has a long history, most notably used by Weierstrass. It is in many ways similar to the modern notion of an algebraic variety, however there is not enough evidence to show that they are indeed the same, so this term has been avoided.

\section*{Acknowledgements.}
\emph{Special thanks go to John Rawnsley, Colin McLarty, John Stillwell, Steve Russ, Derek Holt, Daniel Lewis and Anselm Br\"{u}ndlmayer for their invaluable advice, and in particular to Jeremy Gray for providing this opportunity and for his support throughout.}

\end{document}